\documentclass[12pt]{article}

\usepackage{arxiv}
\usepackage{float}
\usepackage[utf8]{inputenc} 
\usepackage[T1]{fontenc}    
\usepackage{hyperref}       
\usepackage{url}            
\usepackage{booktabs}     
\usepackage{amsmath}
\usepackage{enumitem} 
\usepackage{amsfonts}       
\usepackage{nicefrac}       
\usepackage{microtype}      
\usepackage{lipsum}
\usepackage{graphicx}
\usepackage{algorithm}
\usepackage{algorithmic}
\usepackage{url}

\numberwithin{equation}{section}
\newtheorem{definition}{Definition}[section]
\newtheorem{remark}{Remark}[section]
\newtheorem{assumption}{Assumption}[section]

\graphicspath{ {./images/} }

\title{Electrical Impedance Tomography for Anisotropic Media: a Machine Learning Approach to Classify Inclusions}

\author{
 Romina Gaburro \\
  Department of Mathematics and Statistics\\
  University of Limerick\\
  Health Research Institute (HRI) \\
  \texttt{romina.gaburro@ul.ie} \\
   \And
 Patrick Healy \\
  Department of  Computer Science and Information Systems\\
  University of Limerick\\
  \texttt{patrick.healy@ul.ie} \\
  \And
 Shraddha Naidu \\
 Department of Mathematics and Statistics\\
  University of Limerick\\
  \texttt{shraddha.naidu@ul.ie} \\
  \And
Clifford Nolan \\
  Department of Mathematics and Statistics\\
  University of Limerick\\
  Health Research Institute (HRI) \\
  \texttt{clifford.nolan@ul.ie} \\
}

\begin{document}
\maketitle
\begin{abstract}
We consider the problem in Electrical Impedance Tomography (EIT) of identifying one or multiple inclusions in a background-conducting body
$\Omega\subset\mathbb{R}^2$, from the knowledge of a finite number of electrostatic measurements taken on its boundary $\partial\Omega$ and modelled by the Dirichlet-to-Neumann (D-N) matrix.  Once the presence of one inclusion in $\Omega$ is established, our model, combined with the machine learning techniques of Artificial Neural Networks (ANN) and Support Vector Machines (SVM), may be used to determine the size of the inclusion, the presence of multiple inclusions, and also that of anisotropy within the inclusion(s). Utilising both real and simulated datasets within a 16-electrode setup, we achieve a high rate of inclusion detection and show that two measurements are sufficient to achieve a good level of accuracy when predicting the size of an inclusion. This underscores the substantial potential of integrating machine learning approaches with the more classical analysis of EIT and the inverse inclusion problem to extract critical insights, such as the presence of anisotropy.
\end{abstract}


\section{Introduction}
\label{sec:int}
Electrical Impedance Tomography (EIT) is an affordable and non-invasive imaging modality. Its underlying inverse problem, known as the inverse conductivity problem, or Calder\'on's problem \cite{C}, is about determining the conductivity $\sigma$ of a body $\Omega\subset\mathbb{R}^n$, $n\geq 2$ from the Dirichlet-to-Neumann (D-N) boundary map, which encodes infinite measurements of voltage and electrical current density taken on the boundary $\partial\Omega$ of $\Omega$. In the absence of internal sources or sinks, the electrostatic potential $u$ in a domain $\Omega\subset{\mathbb R}^n$, is governed by the equation
\begin{equation}\label{eq conduttivita'}
   \mbox{div}(\sigma\nabla{u}) = 0 \qquad \mbox{in} \quad \Omega ,
\end{equation}
where $\sigma(x)=(\sigma_{ij}(x)_{i,j=1}^n$, $x\in\Omega$, is a symmetric, positive definite matrix-valued function representing the (possibly anisotropic) electric conductivity of $\Omega$. The inverse conductivity problem consists of finding $\sigma$ when the D-N map $\Lambda_{\sigma}$,
\begin{equation}
\Lambda_{\sigma}:\,{H}^{1/2}(\partial\Omega)\,\ni\,u\vert_{\partial\Omega}\,
\rightarrow\,{\sigma}\nabla{u}\cdot\nu\vert_{\partial\Omega}\in{H}^{-1/2}(\partial\Omega)
\label{eq:ContinuumModel}
\end{equation}
is given for any solution $u \in {H}^{1}(\Omega)$ of (\ref{eq conduttivita'}). Here, $\nu$ denotes the outer unit normal to $\partial\Omega$. This problem arises in many different fields such as geophysics, medical imaging, and non-destructive testing of materials, to mention only a few of its applications. The first mathematical formulation of the inverse conductivity problem is due to Calder\'{o}n \cite{C}, where the problem of uniquely determining the isotropic conductivity $\sigma = \gamma I$ from $\Lambda_\sigma$, was addressed. As main contributions in this respect we mention the papers by Alessandrini \cite{A1}, Kohn and Vogelius \cite{Ko-V1, Koh-V2}, Nachman \cite{N}, and Sylvester and Uhlmann \cite{Sy-U}. We also refer to \cite{Bo, U} for an overview on this inverse problem.

Due to the severe ill-posedness of this inverse problem, even in the isotropic case, EIT is very sensitive to noise perturbations. Alessandrini proved in \cite{A} that, in the isotropic case and dimension $n\geq 3$, assuming \textit{a-priori} bounds on $\sigma$ of the form $\|\sigma\|_{H^s(\Omega)}\leq E$ , $s>(n/2)+2$, leads to a continuous dependence of the estimated $\sigma$ in $\Omega$ upon $\Lambda_{\sigma}$ of logarithmic type. 

Additionally, when the conductivity is anisotropic, its unique determination from the classical D-N map is still an open problem. In fact, it is well-known that for any diffeomorphism $\Phi :\overline\Omega\longrightarrow\overline\Omega$ such that $\Phi\big\vert_{\partial\Omega} = I$, $\sigma$ and its push-forward under $\Phi$,
\begin{equation}
\Phi^{*}\sigma = \frac{(D\Phi )\sigma(D\Phi)^T\big)}{\det(D\Phi)\Big)}\circ \Phi^{-1}
\end{equation}
give rise to the same D-N map \cite{Ko-V1}. Since the above observation, different lines of research have been investigating this inverse problem in the presence of anisotropy. One is that of determining $\sigma$ modulo diffeomorphisms which fix the boundary, in this respect we refer to the seminal paper of Lee and Uhlmann \cite{Le-U} (see also \cite{La-U} \cite{La-U-T} and \cite{Astala2005, Sylvester1990} for the two-dimensional case). 

Anisotropy presents itself in various applications, such as medical imaging or geophysics. In the geophysical context, it was recognised in \cite{S} that anisotropy may affect the electrical properties of geological formations. 
 The loss of anisotropy in tissues is used to detect dystrophic muscle conditions, serving as a potential index of disease status \cite{rutkove2016loss}.  
It is therefore important to investigate under which circumstances $\Lambda_{\sigma}$ uniquely determines an anisotropic $\sigma$. 
 We propose using Artificial Neural Networks (ANNs) to detect the presence of anisotropy within an inclusion. 
 We achieve high accuracy ($>$ 90\%) in detecting special types of anisotropic inclusions in isotropic and anisotropic backgrounds (see Section \ref{sec:AnisoEIT}). 

Detecting an inclusion within a conductor from electrostatic measurements provided by the relevant D-N map is a companion problem to the inverse conductivity problem \cite{AdC}. In 1988 Isakov \cite{I} proved uniqueness for this inverse problem. Regarding the issue of stability, a logarithmic type of stability in terms of the Hausdorff distance between closed sets was given in \cite{AdC} in the case where $\sigma$ is isotropic and homogeneous (constant) within the background and the (possibly disconnected) inclusion. This result was extended in \cite{dC} to the isotropic inhomogeneous conductivity case, where the logarithmic stability of an unknown inclusion with H\"older continuous conductivity (in a smoother isotropic background) is given in terms of a local boundary map.

In the anisotropic case, results of uniqueness and stability for the inclusion problem are quite limited. Logarithmic stability was proved in \cite{dC-R1} for the case of an inclusion $D$ of conductivity $kA$, $x\in D$, where $A$ is a known Lipschitz continuous matrix-valued function representing the background anisotropic conductivity in $\Omega\setminus\overline{D}$ and $k\in \mathbb{R}$ is an unknown positive parameter, which is related to the jump in the conductivity on $\partial D$. Here we train in two dimensions a Radius-Fully Connected Neural Network (R-FCNN) in a similar setup of \cite{dC-R1}, with $A$ a constant diagonal matrix (see \eqref{eq:atank}, \eqref{eq:anisomat}. We then extend this setup to more general settings (see Table \ref{tb:AnisoType}) that include training the network to distinguish between a constant isotropic inclusion $D$ from a smooth spatially varying diagonal anisotropic conductivity (see \eqref{eq:spatialcond}) in an isotropic constant background. With minor modifications, our argument should hold true also when $A$ is merely Lipschitz continuous and non-diagonal. Partial theoretical results in the unique determination of an inclusion and its anisotropic conductivity were also given in \cite{Kw}. 

A further question in this direction is whether the inclusion $D$ can be uniquely (or even stably) determined by a single measurement of voltage and current flux on $\partial\Omega$. Precisely, denoting by $\chi_D$ the characteristic function of $D$, if $u$ is the electrostatic potential in $\Omega$, the problem is to determine $D$ in the equation
\begin{equation}\label{inclusion pb}
\mbox{div}((1+(k-1)\chi_D)\nabla u) = 0\qquad\textnormal{in}\quad\Omega, 
\end{equation}
given one pair (or a finite number) of Cauchy data on $\partial\Omega$,
\begin{equation}\label{Chauchy data}
u\vert_{\partial\Omega}=f;\qquad (\partial u/\partial\nu)\Big\vert_{\partial\Omega} = \psi.
\end{equation}
For the uniqueness issue, partial results have been given in \cite{A2}, \cite{A-I}, \cite{A-I-P}. We also refer to \cite{Alb-San1,Alb-San2} for recent developments in this direction in Calder\'on and related inverse problems. Here we propose the use of Machine Learning (ML) methods, namely Support Vector Machines (SVMs) to detect the presence of an inclusion from the knowledge of the D-N matrix. We are able to consistently detect the presence of inclusion(s) and their size, while estimating their number doesn't seem an easy task (see sections 4-6).

It has been shown that rather than recovering the precise shape of the inclusion $D$, one can stably infer the so-called size of $D$ (its Lebesgue measure) from one measurement $(f,\psi)$ under mild assumptions on the conductivity of both the background $\Omega\setminus\overline{D}$ and the inclusion $D$ (see the seminal work of \cite{Br}, \cite{Fr}). We also refer to \cite{ABRV} (and references therein) for the related problem of (stably) determining an inaccessible part of $\partial\Omega$ in terms of one measurement.  We propose to use ANNs to detect the size of a circular inclusion using the available dataset \cite{data}. The importance of estimating the size of an inclusion is evident in many practical applications. For example, changes in bladder size may indicate fullness and EIT has been successfully applied here \cite{bladderSize}. 

We use a 16-electrode system based on the dataset \cite{data}. Although the theoretical minimum for size determination of an inclusion is one measurement, we investigate whether 16 measurements enhance our ability to detect inclusions when applying ML methods. Given this extra data, we anticipate and find that our results yield a very high level of accuracy in size estimation, further confirming that size can be accurately estimated from the discretised version of the D-N map (the D-N matrix) of EIT. 

We also investigate whether the high level of accuracy in inclusion radii detection can be achieved with fewer measurements from a 16-electrode system. Our findings indicate that a minimum of two measurements is required to determine the inclusion radii in such a system, suggesting that fewer measurements are sufficient for the ANN when the system includes enough electrodes (see Section \ref{sec:IncRad}).

In this proof of concept paper, ML methods are employed to evaluate their performance against traditional analytical techniques for both Calder\`on and the inclusion problem in EIT. The suitability of exploiting the discretised D-N map $\Lambda_{\sigma}$ using ML for a) detecting the presence of a circular inclusion in a 2D body $\Omega$, b) detecting the number of inclusions in $\Omega$, c) detecting the radius of the circular inclusion/s, and d) identifying anisotropy in the inclusions, is explored.  
The experimental setup is chosen to suit the available dataset in \cite{data}, where a circular saline tank, having an inclusion which is \textit{a-priori} known to be circular, is considered. 
The inclusion, assumed to be made of thermoplastic, is constrained to have one of four different radii as in Table \ref{tb:class}.

The paper follows the line of reasoning introduced in \cite{agnelli}, where the D-N matrix was provided as input to an ANN to classify the type of stroke in a brain. Here, the D-N matrix is analysed by various ML methods to determine the presence, number of inclusions, size of inclusion and the presence of anisotropy within the inclusion. We use the continuum electrode model for our analysis (see \eqref{eq:ContinuumModel},\eqref{eq:D-N}). Despite the availability of more intricate electrode models, such as the shut model, gap model, and the complete electrode model, we leverage the computational ease of the continuum electrode model to explore the feasibility and limitations of the method introduced in \cite{agnelli} for the detection of anisotropic inclusions and other inclusion characteristics within a 2D body. Following the argument in \cite{agnelli}, we approximate the D-N map $\Lambda_{\sigma}$ with the discrete D-N matrix $\mathbf{L}_\sigma$. 

For inclusion radii detection, we use the available dataset \cite{data} and adopt the opposite-injection pattern as used in that dataset only. This choice makes the number of data measurements at a manageable level. For anisotropy, due to the lack of available real data at the time of the writing, we simulate the data to match the setup of \cite{data}. 
We adopt the injection method described in \cite{agnelli}, which employs an optimal injection pattern for this situation.
Ideally, we would also use this optimal pattern for our real dataset during inclusion radii detection. However, since the data was collected using the opposite injection pattern, we retain that approach. 

For the case of inclusion size, number and presence detection we use an isotropic tank with an isotropic inclusion. Initially, we do not add anisotropy for simplicity. In the case of anisotropy detection, we constrain the inclusion $D$ to have a special type of anisotropy as in \eqref{eq:aniso}, \eqref{eq:anisomatOffDiag} and \eqref{eq:spatialcond} and employ ANNs to detect it within both an isotropic and an anisotropic background (see Section \ref{sec:AnisoEIT}).

The paper is organized as follows; Section \ref{sec:formulation} provides a detailed formulation of the inverse conductivity problem. Section \ref{sec:ClassModel} describes the classification models used. In Section \ref{sec:incdet} we describe the process of detecting the presence of one and more inclusions (Subsections \ref{sec:1 incl} and \ref{sec:num}, respectively). In Section \ref{sec:IncRad}, we detail the process of using the D-N matrix for detecting inclusion radii. 
We also investigate the impact of reducing the number of measurements/electrodes on inclusion radii detection. Section \ref{sec:AnisoEIT} delves into detecting anisotropy within an inclusion in both an isotropic and anisotropic background. Section \ref{sec:con} concludes our findings and outlines future work. 

\section{Formulation of the problem}
\label{sec:formulation}


\subsection{A-priori information on the domain and the conductivity}\label{subsec:domain and conductivity}
Throughout this paper, we follow the setup of \cite{data}, where the object under investigation $\Omega$ is a 2D disk centred at the origin and with radius 28cm. This setting could easily be adapted to suit a more general domain $\Omega\subset\mathbb{R}^2$ with smooth boundary. In \cite{data}, a 3D tank with the inclusion placed at a fixed height was considered. Here, for the sake of simplicity, we deal directly with a 2D cross-sectional model represented by a disk $\Omega\subset\mathbb{R}^2$. We adopt the following notation:
\begin{enumerate}
\item isotropic conductivities will be denoted by
\begin{equation}
 \sigma^{iso}(x)  = \gamma(x)\mathbf{I},\qquad\textnormal{for\:any}\:x\in\Omega,
\label{eq:iso}
\end{equation}
where $\mathbf{I}\in\mathbb{R}^{2\times 2}$ is the identity matrix and $\gamma$ is a function satisfying, for some constant $\lambda>0$,
\begin{equation}\label{ellipticity isotropic}
\lambda^{-1}\leq \gamma(x)\leq \lambda;\qquad\textnormal{for\:a.e.}\quad x\in\Omega;
\end{equation}
\item anisotropic conductivities will be denoted by
\begin{equation}
\sigma^{aniso}(x) = \mu(x) \sigma(x)
,\qquad\textnormal{for\:any}\:x\in\Omega,
\label{eq:aniso}
\end{equation}
\end{enumerate}
where $\mu$ is a function satisfying \eqref{ellipticity isotropic} and $\sigma$ is an \textit{a-priori} known matrix-valued function, symmetric and satisfying the uniform ellipticity condition
  \begin{equation}\label{ellipticity anisotropic}
    \mathcal{E}^{-1}|\xi|^2 \leq \sum_{i,j=1}^n\sigma(x)\xi_i \xi_j \leq \mathcal{E}|\xi|^2,\qquad\textnormal{for\:a.e.}\quad x\in\Omega,\quad\textnormal{for\:any}\:\xi\in\mathbb{R}^2,
\end{equation}
for some constant $\mathcal{E}>0$.  
 

\subsection{The continuum electrode model and the D-N matrix}\label{subsec:electrodes}
To discretise the D-N map \eqref{eq:D-N}, we consider the continuum electrode model with 16 electrodes, where each electrode positioned along $\partial \Omega$ is of equal area and it is assumed that there are no gaps between any such two electrodes.  

Due to the lack of an available dataset at the time of writing, the data has been simulated via $16$ electrodes for the problems of detecting the presence of an inclusion in $\Omega$; classifying the number of inclusions; and determining the presence of anisotropy within an inclusion. The real data in \cite{data} has been used to detect the size of the inclusion/s, again with 16 electrodes. For the simulated data, the patterns of induced voltage used here are the trigonometric functions in \eqref{eq:vol} as in \cite{agnelli}, where a 32-electrode setup was used. Here, to match the setting adopted in \cite{data} throughout the entire manuscript, a system of 16 electrodes is used for both the real and simulated data sets. For inclusion radii detection, as in \cite{data}, the real data set has been obtained by means of the opposite injection pattern within the 16 electrodes, where current consistently flows between two opposing electrodes  \cite{opp}. 

In \cite{data} as in most practical applications of EIT, currents are injected into the body, and the resulting voltages are measured. In our simulations, we find it more convenient to apply voltages at the boundary using the \texttt{PDE} toolbox in MATLAB2023b and measure the resulting induced currents. 

Pairs of Cauchy data (voltages and surface current densities) are generated through 16 electrodes to form the D-N matrix, as in \eqref{eq:D-N}. In our continuum electrodes model, the trigonometric functions in \eqref{eq:vol} induce voltage patterns simultaneously on all electrodes (see \cite{agnelli}\cite{mueller2012linear}\cite{Syren}). The voltage $V^k$, for  $k \in \{1,2,\dots\}$ is given by
\begin{equation}V^k (\theta) =  \pi^{-1/2}
    \begin{cases}
        \cos\left((k+1)\theta/2\right), & \textnormal{$k$ odd}, \\
        \sin\left(k\theta/2\right), & \textnormal{$k$ even},
    \end{cases}
    \label{eq:vol}
\end{equation}
with $\theta\in[0,2\pi]$. To simulate the boundary current density $J^k$ corresponding to the boundary voltage $V^k$, we solve the forward Dirichlet problem for the interior voltage $u_k$ 
\begin{equation}\label{forward BVP}
    \begin{cases}
        \text{div}(\sigma\nabla u_k) = 0, & \text{in $\Omega$}   \\
        u_k=V^k, & \text{on $ \partial\Omega$},
    \end{cases}
\end{equation}
for $k=1,2,\dots $ and compute the conormal derivative of $u_k$
\begin{equation}\label{J}
J^k\big\vert_{\partial\Omega} := \sigma \nabla u_k\cdot\nu\big\vert_{\partial\Omega},\qquad\textnormal{for}\quad k=1, 2\dots,
\end{equation}
where $\sigma$ is either $\sigma^{iso}$ or $\sigma^{aniso}$ as defined in \eqref{eq:iso} or \eqref{eq:aniso}, respectively. \eqref{forward BVP} is solved using the Finite-Element Method (FEM) in-built in the \texttt{PDE} toolbox in MATLAB2023b, using a discretised FEM mesh of 411 triangles around the domain. The PDE Toolbox in MATLAB also provides the gradient $ \nabla u_k $ for the solution $u_k$. Utilizing the computed $ \nabla u_k$, we can derive $J^k$ as in \eqref{J}.
As the number of electrodes adopted here is $16$, in our continuum electrodes model we have
\begin{equation}\label{regions electrodes}
\partial\Omega = \bigcup_{l=1}^{16} e_l,\quad\textnormal{with}\quad e_l\cap e_m = \emptyset,\quad\textnormal{for}\quad l\neq m.
\end{equation}
Electrode $l$ occupies region $e_l$, and
\begin{equation}\label{measure electrodes area}
|e_l| = |\partial\Omega|/16,
\end{equation}
$l=1,\dots 16$. Here $|e_l|$ and $|\partial\Omega|$ denotes the Lebesgue measure of $e_l$ and $\partial\Omega$, respectively. We discretise $\theta$ to $\theta_l$ as,
\begin{equation}
    \theta_l = 2\pi l/16.
\end{equation}
Here $\theta_l$ denotes the angle the centre of the electrode $l$ makes with the horizontal axis. $V^k$ in \eqref{eq:vol} is discretised into the vector:
\begin{equation}\label{V_discretised_mean} 
\widehat{V}^k = (V^k_l)_{l=1}^{16} \in \mathbb{R}^{16}, \quad \text{with} \quad V^k_l :=  
 V^k(\theta_l),\quad \text{for} \quad k, l = 1, \dots, 16. 
 \end{equation}
Here $\widehat{V}^k_l$ represents the $k^{\text{th}}$ voltage pattern applied to electrode $l$. We also calculate the discretised $J^k$,
\begin{equation}\label{J discretised}
\widehat{J^k} = (J^k_l)_{l=1}^{16} \in \mathbb{R}^{16}, \quad \text{with} \quad J_l^k := (|e_l|)^{-1}\int_{\partial \Omega} J^k, \quad\textnormal{for} \quad k,l=1 ..., 16,
\end{equation}
where $\widehat J^k_l$ is the average (mean) current density measured in the $l^{th}$ region (electrode) due to the voltage pattern $V^k$. The mesh used here provdes $\partial\Omega$ with 400 nodes.
As in \cite{agnelli}, \cite{Syren} we define the $(i,j)^{th}$ entry of the D-N matrix $\mathbf{L}_\sigma\in \mathbb{R}^{16\times16}$ as
\begin{equation}\label{eq:D-N}
    \mathbf{L}_\sigma^{i,j} =  \left< \widehat{J}^j, \widehat{V}^i\right>,
    \end{equation}
for $i,j \in \{1,\dots , 16\}$, where $\left< \cdot, \cdot\right>$ denotes the Euclidean inner product on $\mathbb{R}^{16}$.
As in \cite{agnelli}, we reshape $\mathbf{L_\sigma} \in \mathbb{R}^{16\times16}$ into a column vector $\mathbf{L_\sigma} \in \mathbb{R}^{256}$. In \cite{agnelli}, $\mathbf{L}_{\sigma}$ in \eqref{eq:D-N} is obtained by inversion of the N-D matrix to circumvent the need for numerical differentiation. Here we  find it more convenient to use the \texttt{createPDE} function in MATLAB2023b from the PDE Toolbox for our simulations instead, as this is an in-built MATLAB function that is memory efficient and uses the FEM to solve \eqref{forward BVP}
and form $\mathbf{L}_{\sigma}$ from \eqref{eq:D-N} directly. 

To simulate noise in the data collection process of our simulated data, 
 we add a zero-mean Gaussian noise to $\mathbf{L}_\sigma$ 
\begin{equation}\label{eq:ND-N}
    \mathbf{L}^{'}_\sigma = \mathbf{L}_\sigma + 10^{-2}\:\mathbf{N}.
    \end{equation}
The noise vector $\mathbf{N} \in \mathbb{R}^{256}$ is generated using the \texttt{randn()} function in MATLAB2023b, so that $\mathbf{N} \sim \mathcal{N}(0,1)\,$. We scale the noise $\mathbf{N}$ in \eqref{eq:ND-N} by $10^{-2}$ as a noise level of $1\%$ is typically used in industrial EIT \cite{noise}. We do not add noise to the measurements from the real dataset \cite{data} used for the radii-detection problem in Section \ref{sec:IncRad} as this data was collected using an EIT system which already included noise in the measurements. Additionally, in the real data, due to the nature of the measurements taken, we calculate the D-N matrix as $\mathbf{L}_\sigma = \mathbf{R}^{-1}_\sigma$, where $\mathbf{R}_\sigma$ is the Neumann-to-Dirichlet (ND) matrix. Next, $\mathbf{L}^{'}_\sigma$ and $\mathbf{L}_\sigma$ in \eqref{eq:ND-N} and \eqref{eq:D-N}, respectively, are used as input to our classification ML models (Artificial Neural Networks and Support Vector Machines) to detect the presence of inclusions, the number of multiple inclusions, the size of such inclusions and the presence of anisotropy within them.  We first give a brief introduction to the ML classification models considered here. The pseudocode algorithm for D-N matrix generation is presented below in Algorithm \ref{alg:d-n-gen}.

\begin{algorithm}[H]
\caption{D-N matrix and dataset simulation}
\begin{algorithmic}[1]
\FOR{each sample in dataset}
\STATE Simulate a 2D tank domain with/without inclusion using \texttt{domain=decsg()} in MATLAB for class label $N$.
\STATE Use \texttt{model = createpde()} in MATLAB to generate a PDE object.
\STATE Define $\sigma$ for the tank and/or inclusion(s).
\STATE Generate FEM Mesh with 411 triangular elements for the tank/inclusion domain using \texttt{geometryFromEdges(model,domain)} and \texttt{generateMesh(model)}.
\FOR{$k = 1 \text{ to } 16$}
    \STATE Apply the $k^{\text{th}}$ voltage pattern \eqref{eq:vol} as the boundary condition.
    \STATE Find the voltage at the center of each electrode $\hat{V}^k_l:=   V^k(\theta_l)$.
    \STATE Solve the PDE using FEM with \texttt{solvepde(model)}
    \STATE Store  results: $u_k$ and $\nabla u_k$.
    \STATE Calculate the resulting current density $J^k=\sigma \nabla u_k\cdot\nu\big\vert_{\partial\Omega}$.
    \STATE Find the mean current density at each electrode:  \\$\widehat{J^k} = (J^k_l)_{l=1}^{16} \quad, J_l^k := (|e_l|)^{-1}\int_{\partial \Omega} J^k$.
\ENDFOR
\STATE Calculate D-N matrix $\mathbf{L}_\sigma(i,j) =  \left< \widehat{J}^j, \widehat{V}^i\right>$.
\STATE Reshape the noisy D-N matrix $\mathbf{L}_\sigma$ into a row vector to use as input for classification.
\STATE Add noise to the D-N matrix $\mathbf{L}^{'}_\sigma = \mathbf{L}_\sigma + 10^{-2}\cdot \mathbf{N}.$
\STATE Store $\{\mathbf{L}^{'}_\sigma, n \}$ as data-point.
\ENDFOR
\STATE Train ML algorithm using data: $\{\mathbf{L}^{'}_\sigma, n \}$.
\STATE Classify $\mathbf{L'_{\sigma}} \quad\xrightarrow[]{\text{ML Algorithm}} n$.
\end{algorithmic} \label{alg:d-n-gen}
\end{algorithm}


\section{The ML classification models adopted here}
\label{sec:ClassModel}


\subsection{Artificial Neural Networks}
\label{subsec:ann}

We employ an Artificial Neural Network (ANN) with one hidden layer to classify the D-N matrix (obtained either via simulated or real data as per algorithms \ref{alg:d-n-gen} or \ref{alg:d-n-real}, respetively) corresponding to the experiments described in Section \ref{sec:formulation}:
\begin{equation}
    \text{ANN: D-N Matrix} \mapsto \mathrm{Classification}.
\end{equation}
To introduce non-linearity in the network, an activation function, such as the \textit{sigmoid function} $\phi: \mathbb{R}\longrightarrow [0,1]$, defined by 
\begin{equation}
\phi(z) = \left(1 + e^{-z}\right)^{-1},
\label{eq:sig}
\end{equation}
is applied. Here our neural network is represented by a function
\begin{equation}\label{ANN function}
    \mathcal{F}: \mathbb{R}^{256}\longrightarrow\mathbb{R},
\end{equation} 
having an input layer with 256 neurons - each input neuron corresponding to a component of either $\mathbf{L_\sigma} \in \mathbb{R}^{256}$ or $\mathbf{L'_\sigma} \in \mathbb{R}^{256}$, one hidden layer with 64 neurons and an output layer with $n$ neurons, where $n$ is the number of output possible classifications of the data $\mathbf{L_\sigma}$ or $\mathbf{L'_\sigma}$. A probability function, known as the \textit{ softmax function} \eqref{eq:soft}, is applied component wise to the final layer of $n$ neurons, allowing for the selection of the neuron in the layer with the highest probability, resulting in the prediction $O_{net}$, a label used to classify our output according to Definition \ref{def:Onet} below. 
If we denote each input vector $\mathbf{L_\sigma}^T$ or $\mathbf{L'_\sigma}^T$ by $\mathbf{x}^{(1)} \in \mathbb{R}^{1\times 256}$, the matrix of weights on the first weight layer by $\mathbf{W}^{(1)} \in \mathbb{R}^{256\times64}$, the vector of biases on the hidden layer by $\mathbf{B}^{(1)} \in \mathbb{R}^{1\times 64}$ (as there are 64 neurons in the hidden layer), then the input to the second layer is given by 
\begin{equation}
    \mathbf{x}^{(2)} = \phi(\mathbf{x}^{(1)}\mathbf{W}^{(1)} + \mathbf{B}^{(1)}).
    \label{eq:input_nn}
\end{equation}
In \eqref{eq:input_nn} it is understood that, with abuse of notation, $\phi$, defined in \eqref{eq:sig}, is applied component wise to $\mathbf{x}^{(1)}\mathbf{W}^{(1)}+ \mathbf{B}^{(1)}\in\mathbb{R}^{1\times 64}$, leading to $ \mathbf{x}^{(2)}\in \mathbb{R}^{1\times 64}$. Let $\mathbf{W}^{(2)} \in \mathbb{R}^{64\times n}$ denote the weights connecting the hidden layer to the output layer (with 64 neurons in the hidden layer and $n$ neurons in the output layer) and $\mathbf{B}^{(2)} \in \mathbb{R}^{1 \times n}$, leading to $\mathbf{x}^{(3)} \in \mathbb{R}^{1 \times n}$. Here, $n$ is the number of possible classifications of the data and $O_{net}$ introduced above can be any of the $n$ components of $\mathbf{x}^{(3)}$. Subsequently, we apply the softmax function to $\mathbf{x}^{(3)} = ( x^{(3)}_i)_{i=1}^n$ component wise as,
\begin{equation}
     \varphi(x^{(3)}_i) = (e^{x^{(3)}_i})/(\sum_{j=1}^n e^{x^{(3)}_i}),\qquad i={1,..,n},
     \label{eq:soft}
 \end{equation}
to obtain the probability distribution indicating the likelihood of input $\mathbf{L}^{'}_\sigma / \mathbf{L}_\sigma$ to belong to one of the $n$ \textit{classes}. 
\begin{definition}\label{def:Onet}
We will say that the input of the network belongs to the class $J$, for some $1\leq J \leq n$ if 
\begin{equation}
    \max_{i=1,\dots , n} \varphi (x^{(3)}_i) = \varphi (x^{(3)}_J).
\end{equation}
If that is the case we will denote the corresponding output with
\begin{equation}
O_{net}:=J.
\end{equation}
In the case where $J$ is not unique for $x^{(3)}$, one $J$ is randomly selected by the ANN algorithm as the output.
\end{definition}
The classes used in our classifications are summarised in Table \ref{tb:allclass}.
\begin{table}[H]
\centering
\begin{tabular}{llc}
\textbf{Task}                                      & \textbf{Classifications}& \textbf{Class $J$} \\ \hline
\multicolumn{1}{l|}{Inclusion presence detection}  & No Inclusion       & 1                      \\
\multicolumn{1}{l|}{}                              & One Inclusion      & 2                       \\ \hline
\multicolumn{1}{l|}{Number of inclusions detection} & One inclusion      & 1                      \\
\multicolumn{1}{l|}{}                              & Two inclusions     & 2                       \\
\multicolumn{1}{l|}{}                              & Three inclusions   & 3                       \\ \hline
\multicolumn{1}{l|}{Inclusion radii detection}     & Diameter = 19.4 mm & 1                       \\
\multicolumn{1}{l|}{}                              & Diameter = 38.8mm  & 2                       \\
\multicolumn{1}{l|}{}                              & Diameter = 58.2mm  & 3                       \\
\multicolumn{1}{l|}{}                              & Diameter = 77.6mm  & 4                       \\ \hline
\multicolumn{1}{l|}{Anisotropy detection}          & Isotropic          & 1                       \\
\multicolumn{1}{l|}{}                              & Anisotropic        & 2                      
\end{tabular}
\caption{Table showing the classifications used here along with the corresponding class labels ($J$).}
\label{tb:allclass}
\end{table}
Our network is trained
\begin{equation}
    \mathcal{F}(\mathbf{W},\mathbf{B}): \mathbf{x^{(1)}} \longmapsto O_{net},
\end{equation}
to find its optimal parameters (weights ($\textbf{W}$) and biases ($\textbf{B}$)), by minimizing the cross-entropy loss function,
\begin{equation}
    L(\mathcal{F}(\mathbf{W},\mathbf{B})) = - N^{-1} \sum_{i=1}^N \sum_{k=1}^n y_{i, k} \log(\hat{y}_{i, k}).
    \label{eq:loss}
\end{equation}
Here, $n$ and $N$ are the total number of classes and training samples, respectively. $y_{i,k}$ denotes the true label indicator during training $i$, for $i=1,\dots , N$: is 1 if the input $\mathbf{x}^{(1)}$ belongs to class $k$ and 0 otherwise. $\hat{y}_{i,k}$ denotes the predicted likelihood during training $i$, for $i=1,\dots , N$, that $\mathbf{x}^{(1)}$ belongs to class $k$: it is calculated by the model using the softmax function in \eqref{eq:soft}. 

Weights ($\mathbf{W}$) and biases ($\mathbf{B}$) are adjusted iteratively via the Scaled Conjugate Descent (SCD), based on the gradient of the error/loss function with respect to the weights and biases respectively, to solve
\begin{equation}\label{min 1}
    (\mathbf{W}^*, \mathbf{B}^*) = \arg \min_{\mathbf{W}, \mathbf{B}} L(\mathcal{F}(\mathbf{W},\mathbf{B})).
\end{equation}
The choice of the Scaled Conjugate Descent as the minimisation algorithm in \eqref{min 1} is dictated by its computational efficiency in classification tasks.
\begin{definition}
 We will say that the ANN is \textit{trained} when the minimisation algorithm converges and $(\mathbf{W}^*, \mathbf{B}^*)$ is found.
   \end{definition}
\begin{remark}\label{datasplitANN}
For ANNs, we consistently use 80\% of the our datasets for training the network, referred to as the \textit{training data}. The remaining 20\% is split equally: 10\% is used for validation (\textit{validation dataset}) and 10\% for testing (\textit{testing dataset}) (see Table \ref{tb:split}). 
\end{remark}
The validation dataset helps estimate the model's performance during training while adjusting weights and biases. The test dataset, withheld during training, provides an unbiased evaluation of the final model's performance, allowing for comparison and selection among different trained networks. By evaluating the network on this unseen data, we can estimate its true accuracy and gauge how well the model performs on real-world, out-of-sample data. A brief description of the ANN algorithm is described below in Algorithm \ref{ANN algorithm}.
\begin{algorithm}[H]
\caption{ANN}
\begin{algorithmic}[1]
\STATE Initialize $\mathbf{W^{(i)}}$, $\mathbf{B^{(i)}}$
\STATE $index = 1$
\WHILE{$index \leq 1000$}
    \FOR{$i = 1$ \TO $N$} 
        \STATE $\mathbf{x}^{(1)} = \mathbf{L}_{\sigma^i}$
        \STATE $\mathbf{x}^{(2)} = \varphi(\mathbf{x}^{(1)} \mathbf{W^{(1)}} + \mathbf{B^{(1)}})$
        \STATE $\mathbf{x}^{(3)} = \varphi(\mathbf{x}^{(2)})$
        \STATE $\hat{y}_{i,k} = \varphi(\mathbf{x}^{(3)}_k)$
        \STATE Calculate $L$ using $\hat{y}_{i,k},\mathbf{W},\mathbf{B}$
        \STATE Update $\mathbf{W}$, $\mathbf{B}$ using SCD
    \ENDFOR
    \STATE $index = index + 1$
    \IF{L starts increasing per epoch}
        \STATE \textbf{Stop}
    \ENDIF
\ENDWHILE
\end{algorithmic}\label{ANN algorithm}
\end{algorithm}
\begin{figure}[H]
    \centering
    \includegraphics[width=0.43\linewidth,keepaspectratio]{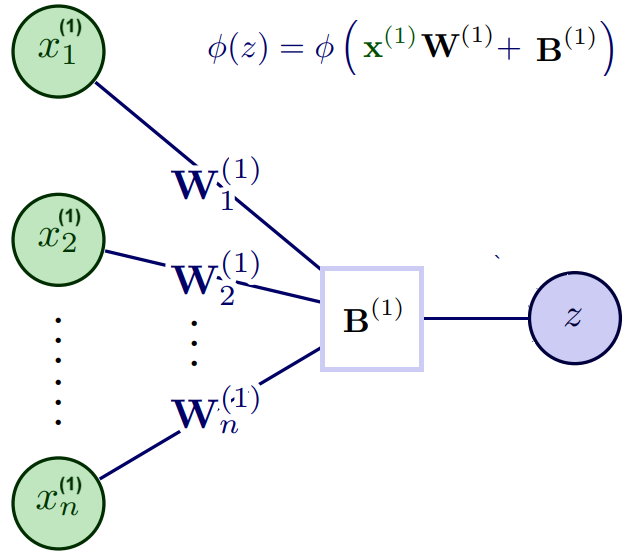}
    \caption{Illustration of the function of a single neuron $z$ (highlighted in a blue circle)  in a neural network. The green circles represent neurons from the previous (input) layer, and the lines indicate the weights and biases connecting each input neuron to the current neuron $z$. Neuron $z$ receives the input $\mathbf{x}^{(1)} = (x_i^{(1)})_{i=1}^n$ (as in \eqref{eq:input_nn}) from the previous layer, where each $\textbf{W}_i^{(1)}$ denotes column $i$ of matrix $\textbf{W}^{(1)}\in\mathbb{R}^{n\times m}$ ($\textbf{W}_i^{(1)}$ corresponds to the weight associated with input $x_i^{(1)}$). $\textbf{B}^{(1)}$represents the bias associated with the current neurons in $\mathbf{x}^{(1)}$, contributing to the neuron output, $z$. The activation function $\phi$ is then applied to $z$ as in \eqref{eq:input_nn}, producing the output $\phi(z)$. If there are subsequent layers, $\phi(z)$ acts as the input $\mathbf{x}^{(2)}$ to the next layer, with an associated weight.} 
   \label{fig:an}
\end{figure}


\subsection{Support Vector Machines}
\label{subsec:SVM_alg}
To detect the presence and number of inclusions, with data given by $\mathbf{L'_{\sigma}} \in \mathbb{R}^{256}$ in \eqref{eq:ND-N}, we use Support Vector Machines (SVMs). The goal of our SVM is to find a \textit{hyperplane} that divides the data points into distinct regions of the \emph{feature space}, $\mathbb{R}^{\ell}$, where $\ell=256$ in this instance. 

\begin{figure}[H]
    \centering
\includegraphics[width=0.8\linewidth]{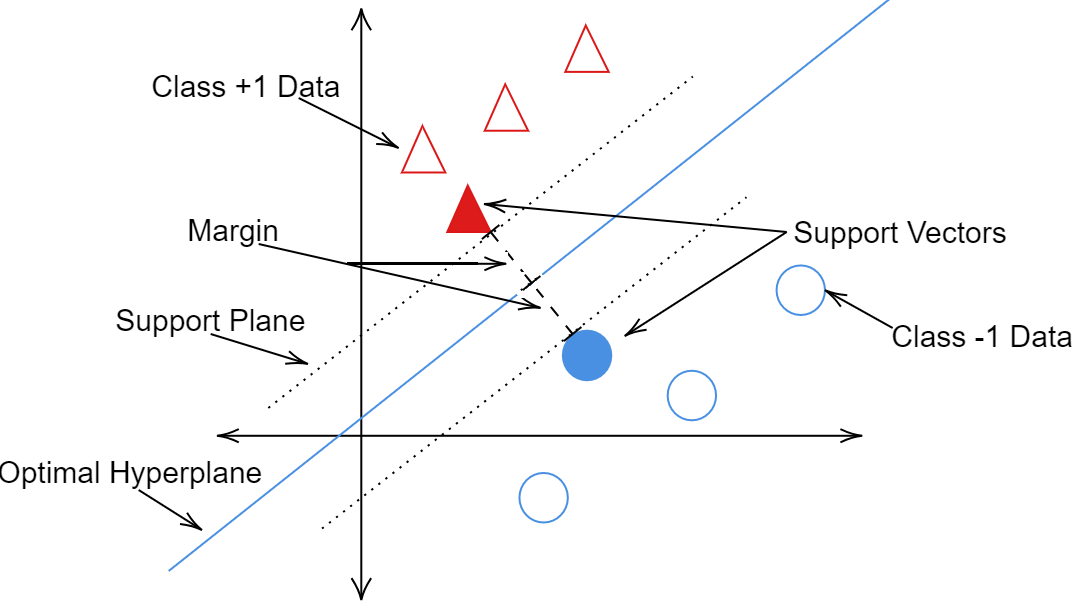}
    \caption{Simplified illustration of the construction of the optimal hyperplane (line) in SVM. In $2D$/$3D$ the \textit{support lines/planes} passing through the support vectors and parallel to the optimal hyperplane provide a simple illustration of the regions of classification and the margin in the SVM classification process.}
    \label{fig:SVMDef}
\end{figure}

For $\mathbf{L'_{\sigma^i}}\in \mathbb{R}^{256}$, $i=1,\dots, N$, set of data points, their classification $h(\mathbf{L'_{\sigma^i}})$, is defined by 

\begin{equation}
h(\mathbf{L'_{\sigma^i}}) = \begin{cases}
    +1,  & \mathbf{w}\cdot\mathbf{L'_{\sigma^i}}  +b \geq +1, \\
    -1, & \mathbf{w}\cdot\mathbf{L'_{\sigma^i}}  +b \leq -1, 
\end{cases}  
\label{eq:svm_class}
\end{equation}
for $i=1,\dots, N$, where $\mathbf{w}$ is the weights vector and $b$ is the bias term. The values of $\mathbf{w}$ and $b$ are determined using \eqref{eq:lag} below. Here a point above the hyperplane will be classified as +1 and below the hyperplane will be classified as -1. Therefore our SVMs' training data for $N$ data points takes the form,
\begin{equation}
    \{\mathbf{L'_{\sigma^i}},y_i\}, \quad \textnormal{where} \quad y_i \in \{-1,+1\}, \quad i\in\{1,2,...,N\},
\end{equation}
where $y_i$ is the output label associated with $\mathbf{L_{\sigma^i}}$. Here, the SVM aims to find the \textit{optimal hyperplane} separating the two classes. 
We give a brief description of the algorithm's steps for the sake of completeness:
\begin{enumerate}
    \item identify the \textit{support vectors}: these are the data points from opposing classes with minimum Euclidean distance;
 \item identify the \textit{optimal hyperplane} by solving the $l^2$ constraint optimisation problem for the weights $\mathbf{w}$ 
\begin{equation}
    M\!\!=\!\min\!\left(\|\textbf{w}\|^2/2\right), \quad \textnormal{such\: that} \quad y_i (\mathbf{w}\cdot\mathbf{L'_{\sigma^i}}+b)\geq 1, \quad\textnormal{for}\quad i= 1,2,...,N,
    \label{eq:optimization}
\end{equation}
to maximize the so-called \textit{margin}, the distance between the support vectors and the optimal hyperplane. A simplified illustration is given in Figure \ref{fig:SVMDef}. The solution to \eqref{eq:optimization} is obtained by solving 
\begin{equation}
    \min_{\mathbf{w},b} \max_{\alpha_i} \mathcal{L}(\mathbf{w},b,\alpha),
    \label{eq:lag1}
\end{equation}
where $\mathcal{L}(\mathbf{w},b,\alpha)$ is the Lagrangian
\begin{equation}
    \mathcal{L}(\mathbf{w},b,\alpha) = \left(\|\textbf{w}\|^2/2\right) - \sum_{i=1}^n \alpha_i [y_i (\mathbf{w}\cdot\mathbf{x}^i +b)],
    \label{eq:lag}
\end{equation}
where $\alpha = (\alpha_i)_{i=1}^n$, with $\alpha_i \geq 0$, are the Lagrange multipliers. 
\eqref{eq:lag1}-\eqref{eq:lag} can be solved via either Quadratic Programming or Sequential Minimal Optimisation (SMO). As we use MATLAB2023b, the default SMO optimiser is used here (see \cite{svm_opti}). 
\end{enumerate}
\begin{remark}\label{def:multiSVM}
In multi-class classification, a hyperplane is created between each pair of classes. As a result, the SVM functions as a recursive binary classifier. For $n$ classes, there are $n(n-1)/2$ hyperplanes. This approach is known as the One-Vs-One SVM. During classification, a data point is evaluated between each pair of classes and the final classification is determined by the majority vote across all pairs.
\end{remark}
We employ a \textit{Quadratic Kernel} for our SVM, as these kernels are well-known for enhancing SVM performance \cite{svmperf} as they capture the interactions between different components of the input vector through \textit{interaction terms}. These terms reveal how one component influences another, enabling the identification of latent patterns and inter-dependencies within the data, $\mathbf{x} \in \mathbb{R}^{\ell}$, by embedding the data into a higher dimensional space. The quadratic kernel $K$ for the SVM is,
\begin{equation}\label{K}
  K(\mathbf{L'_{\sigma^i}},\mathbf{L'_{\sigma^j}}) =\Big((\Phi(\mathbf{L'_{\sigma^i}}))^T \Phi(\mathbf{L'_{\sigma^j}})+1\Big)^2,\quad\textnormal{for}, \quad i,j=1,\dots , N,
\end{equation}
where,
$$\Phi:\mathbb{R}^{\ell} \longmapsto \mathbb{R}^{2\ell+ (\ell(\ell-1))/2}.$$  As $\ell=256$ here, we have $\Phi:\mathbb{R}^{256} \longmapsto \mathbb{R}^{33,152}$, defined by
\begin{eqnarray}
& & \hspace{-0.5cm}\Phi(\mathbf{x})\\
& & \hspace{-0.5cm} = \left[x_1,...,x_{256},x_1^2,\sqrt2x_1x_2,...,\sqrt2x_1x_{256}, x_2^2,\sqrt2x_2x_3,...,\sqrt2 x_2x_{256},....., x_{256}^2\right],\nonumber
\end{eqnarray}
for any $\mathbf{x}\in\mathbb{R}^{256}$, with $\mathbf{x} =(x_k)_{k=1}^{256}$. 
The interaction terms $x_p x_q , p \neq q $ capture the degree to which different components of the input vector $\mathbf{x}$ influence each other. This is particularly useful in problems where the combined effect of two features might reveal important patterns or interactions that aren't obvious from the features alone. Each interaction term appears only once as a component of $\Phi(\mathbf{x})$ in the higher kernel space and is rescaled by $\sqrt{2}$ to prevent it from becoming disproportionately large compared to the individual $x_i$ terms. This ensures that the interaction terms are not ignored during training, as they might otherwise contribute less due to their smaller magnitude when multiplied by the weights. Using $K$ in \eqref{K}, $\mathbf{L_{\sigma}'}$ is embedded onto a higher dimensional space before being classified via the SVM in \eqref{eq:svm_class} through the optimal hyperplane in $\mathbb{R}^{33,151}$ determined by \eqref{eq:lag1}. 
The SVM algorithm is written in pseudocode as:
\begin{algorithm}[H]
\caption{SVM}
\begin{algorithmic}[1]
\STATE Identify the closest points from opposing classes as support vectors using Euclidean distance.
\STATE Embed the input data to the kernel space using: 
\\$
K(\mathbf{L'_{\sigma^i}}, \mathbf{L'_{\sigma^j}}) = \left( \left( \Phi(\mathbf{L'_{\sigma^i}}) \right)^T \Phi(\mathbf{L'_{\sigma^j}}) + 1 \right)^2, \quad \text{for} \quad i,j = 1, \dots, N
$.
\STATE Set the $l^2$ optimization problem to maximize the margin as:\\
$
M = \min \left(\|\textbf{w}\|^2/2 \right), \quad \text{such that} \quad y_i \left( \mathbf{w} \cdot \mathbf{L'_{\sigma^i}} + b \right) \geq 1, \quad \text{for} \quad i = 1: N$.
\STATE Solve the optimization problem using the Lagrangian and SMO:
$
\min_{\mathbf{w}, b} \max_{\alpha_i} \mathcal{L}(\mathbf{w}, b, \alpha)
$
\STATE Draw the hyperplane in $\mathbb{R}^{33,151}$.
\STATE Classify points based on the hyperplane.
\end{algorithmic}
\end{algorithm}
\begin{remark}\label{datasplitSVM}
For SVMs, we consistently use 90\% of the data for training and 10\% for testing. For validation, we applied 5-fold validation as per Definition \ref{K fold} below (see Table \ref{tb:split}). 
\end{remark}
\begin{definition}\label{K fold}
    In $K$-fold validation, the training dataset is split into $ K $ subsets (or folds) of equal size. In each iteration, one of the subsets is used for validation while the remaining $ K-1$ subsets are used for training. 
\end{definition}


\subsection{Performance of our classification models}\label{subsec: performance}
To evaluate and compare the performance of our classification models (ANNs and SVMs) we use a \textit{confusion matrix} (CF), a widely-used tool for assessing classification accuracy. A CF is a square matrix $\mathbf{C}$ of size $n \times n$, where $n$ is the number of classes adopted in the classification task. Each entry $C_{ij}$ represents the number of times that samples from the actual class $i$ were predicted as class $j$. The rows of $\mathbf{C}$ correspond to the actual classes (the true class labels of the samples), where the columns correspond to the predicted classes (the labels assigned by the model). 

The confusion matrices for our trained models on the respective Test Datasets are shown in Figures \ref{fig:QuadSVM}, \ref{fig:QuadSVM_NumInc}, \ref{fig:NumInc388}, \ref{fig:TestC},  \ref{fig:anisodet} and \ref{fig:IsoTankANisoIso}.
The training, testing and validation splits for our models are given in Table \ref{tb:split}, showing that 10\% of the dataset is consistently used for testing.

\begin{table}[htb]
\centering
\begin{tabular}{llll}
\textbf{Model} & \textbf{Training Dataset}& \textbf{Testing Dataset}& \textbf{Validation Dataset}\\ \hline
ANN            & 80\%              & 10\%             & 10\%                \\
SVM            & 90\%              & 10\%             & $5$-Fold         
\end{tabular}
\caption{Training, Testing and Validation split for datasets used in ANNs and SVMs.}
\label{tb:split}
\end{table}


\section{Inclusion presence and inclusions' number detection}\label{sec:incdet}
Through the entire manuscript, the object under investigation (see Section \ref{sec:formulation}) $\Omega$ is a $2D$ disk centered at the origin and with radius 28cm. We will refer to $\Omega$ as the tank. We make the following main assumptions.
\begin{assumption}\label{assumption 1}
We assume that
\begin{enumerate}[label=(\arabic*)] 

\item The conductivity in $\Omega$ is isotropic and homogeneous (constant) with conductivity denoted by
\begin{equation}
 \sigma^{iso}_{tank}  = \lambda_{tank}
\begin{bmatrix}
1 & 0 \\
0 & 1  \\
\end{bmatrix},
\label{eq:itank}
\end{equation}
with $\lambda_{tank}=1.45 \,\text{S} \,\text{m}^{-1}$ \cite{saline}. 
\item If there is an inclusion in $\Omega$, then the conductivity in the inclusion is isotropic and homogeneous (constant) and we denote it by 
\begin{equation}
 \sigma^{iso}_{inc}  = \lambda_{inc}
\begin{bmatrix}
1 & 0 \\
0 & 1  \\
\end{bmatrix}, 
\label{eq:iinc}
\end{equation}
where $\lambda_{inc}\neq\lambda_{tank}$ and, based on typical values for thermoplastics, is assumed to be a constant and to satisfy
\begin{equation}\label{mucinc}
8\: \text{S} \,\text{m}^{-1}\leq\lambda_{inc}\leq 10\:\text{S} \,\text{m}^{-1}
\end{equation}
as in \cite{saline}.

\item If there is an inclusion in $\Omega$, then it is circular and its diameter can assume one of the four values listed in Table \ref{tb:class}, as in \cite{data}
\begin{table}[H]
\centering
\begin{tabular}{c|l}
 \hline
\textbf{Class Label}  & \textbf{Diameter} \\ \hline
1 & 19.4mm            \\ 
2 & 38.8mm              \\ 
3 & 58.2mm              \\ 
4 & 77.6mm              \\ \hline
\end{tabular}
\caption{List of the possible diameters of the possible inclusion and the corresponding class labels for inclusion detection dataset.}
\label{tb:class}
\end{table}
\end{enumerate}
\end{assumption}


\subsection{Detection of a single inclusion}\label{sec:1 incl}
Here we further assume that:
\begin{assumption}\label{assumption 2}
If $\Omega$ has an inclusion, then the number of inclusions in $\Omega$ is at most $1$.
\end{assumption}
To detect whether the tank is empty or contains a single inclusion, we simulate 4800 tank samples in total, of which 2400 contain a circular inclusion and the remaining samples contain no inclusion. Data are simulated using the continuum model \eqref{eq:ContinuumModel} with $16$ electrodes as described in Section \ref{subsec:electrodes}. 
Voltages \eqref{eq:vol} are prescribed at the boundary $\partial \Omega$, and the resulting current densities at $\partial\Omega$ are calculated to form the D-N matrix $\mathbf{L_\sigma}$ in \eqref{eq:D-N}. Noise is added to $\mathbf{L_\sigma}$ according to \eqref{eq:ND-N} to form $\mathbf{L'_\sigma}$ for each simulation. 
Using the quadratic SVM as described in Section \ref{subsec:SVM_alg}, these noisy matrices $\mathbf{L'_\sigma}$ are classified as:
\begin{equation}
\mathbf{L'_{\sigma}} \quad\xrightarrow[]{SVM} \quad
    \begin{cases}
    \text{No inclusion} &  \text{Class 1}\\
    \text{One inclusion} & \text{Class 2}
\end{cases}.
\end{equation}
\begin{figure}[H]
    \centering
    \includegraphics[width=0.5\linewidth]{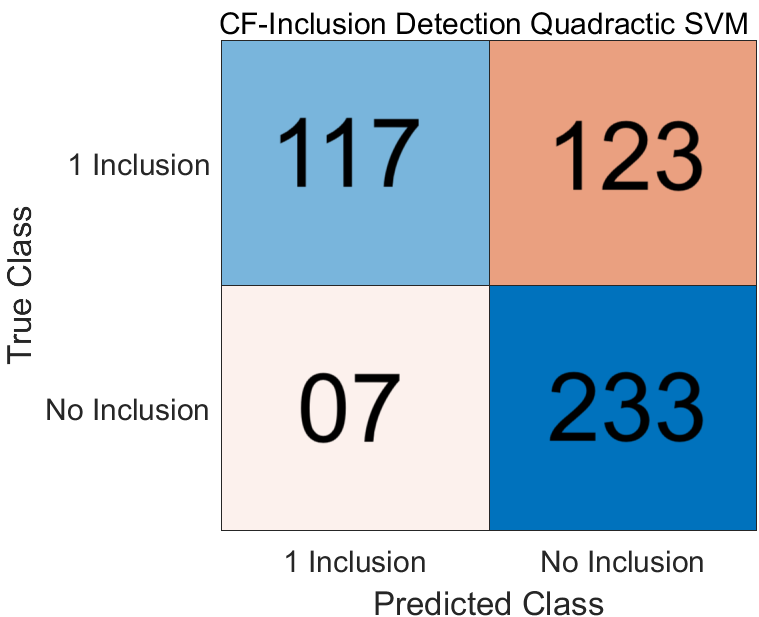}
    \caption{Confusion matrix for Quadratic Kernel, SVM, for a single inclusion detection showing 72.9\% accuracy in prediction. Correct classifications are shown in blue and light blue along the main diagonal of the confusion matrix, while incorrect ones are displayed in red in the off-diagonal entries of the matrix. The darker the colour (blue if the classification is correct), the higher the number of classifications in that category. In an ideal confusion matrix, we would see a dark diagonal (indicating high accuracy) and light red off-diagonal elements (indicating fewer misclassifications).}
    \label{fig:QuadSVM}
\end{figure}
Training the Quadratic Kernel SVM on this data using \texttt{fitcsvm()} on MATLAB2023b, a prediction accuracy of 69.7\% is achieved on the validation dataset. Here, a 5-fold validation according to Definition \ref{K fold}, has been used. When tested on the held-out test dataset, the SVM achieves a final accuracy of 72.9\%. The confusion matrix for the test dataset is given in Figure \ref{fig:QuadSVM}.

Next, we consider the problem of detecting the correct number of inclusions in $\Omega$, when $\Omega$ is allowed to contain more than one inclusions.


\subsection{Classification of the number of  inclusions}
\label{sec:num}
  
Here we make the following main assumption:
\begin{assumption}\label{assumption 3}
There is at least one inclusion in $\Omega$.
\end{assumption}
To evaluate whether the noisy D-N matrix $ \mathbf{L}^{'}_\sigma$ in \eqref{eq:ND-N} contains the necessary information for an ML algorithm like SVM to accurately classify the number of inclusions present in $\Omega$, we simulate the data as in Section \ref{sec:1 incl} via the continuum electrode model. Noise is added to $\mathbf{L_\sigma}$ to form the noisy D-N Matrix $\mathbf{L'_\sigma}$ as in \eqref{eq:ND-N}, which is then fed into the quadratic SVM of Section \ref{subsec:SVM_alg}. We ran two experiments, based on the size the inclusion/s are allowed to have. In the first experiment, the inclusions are kept at a fixed radius $19.4$ mm (Class 1 in Table \ref{tb:class}), opting for the class of smallest inclusions to allow for multiple placements within the tank while maintaining adequate spacing. Three different scenarios are simulated to reflect the presence of 1, 2, or 3 inclusions in the tank leading to the classification problem
\begin{equation}
\mathbf{L'_{\sigma}} \quad\xrightarrow[]{SVM} \quad
    \begin{cases}
    \text{1 Inclusion} &   \text{Class 1}\\
    \text{2 Inclusions} & \text{Class 2}\\
       \text{3 Inclusions} & \text{Class 3}
\end{cases}.
\label{eq:NumInc}
\end{equation}
For multiple classifications, SVM of Section \ref{subsec:SVM_alg} draws hyperplanes between each pair of classes and acts as a recursive binary classifier - see Remark \ref{def:multiSVM}. Due to the presence here of more classes (compared to the single inclusion detection problem of Section \ref{sec:1 incl} with $2$ classes only), more data is generated in this experiment. A total of $6000$ D-N matrices are simulated: $2000$ data points for each of the three scenarios presented in Table \ref{tb:class}.  

\begin{figure}[H]
    \centering
    \includegraphics[width=0.6\linewidth]{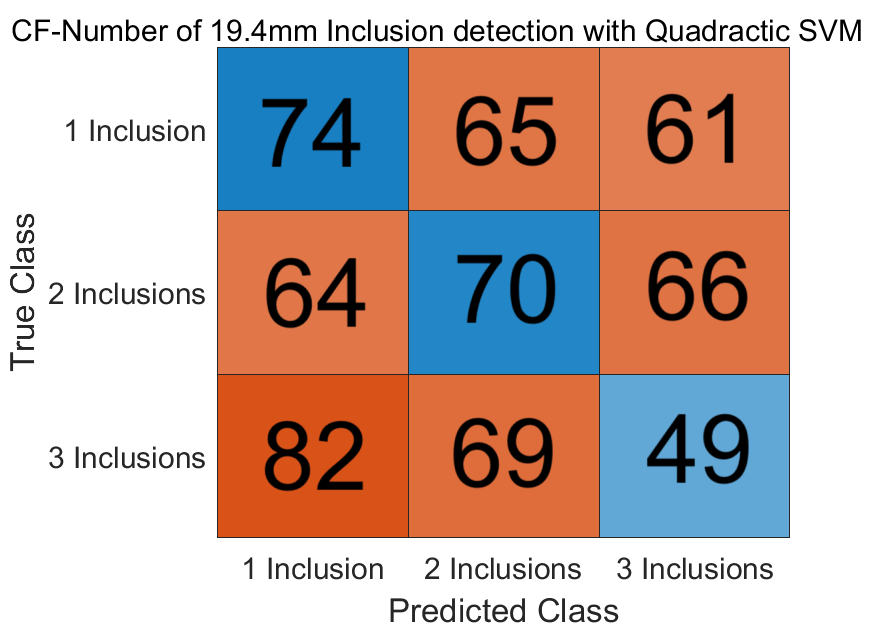}
    \caption{Confusion matrix for the test data for detecting the number of $19.4$mm inclusions in a tank showing $32.2\%$ accuracy in prediction.}
    \label{fig:QuadSVM_NumInc}
\end{figure}

\begin{figure}[H]
    \centering
    \includegraphics[width=0.6\linewidth]{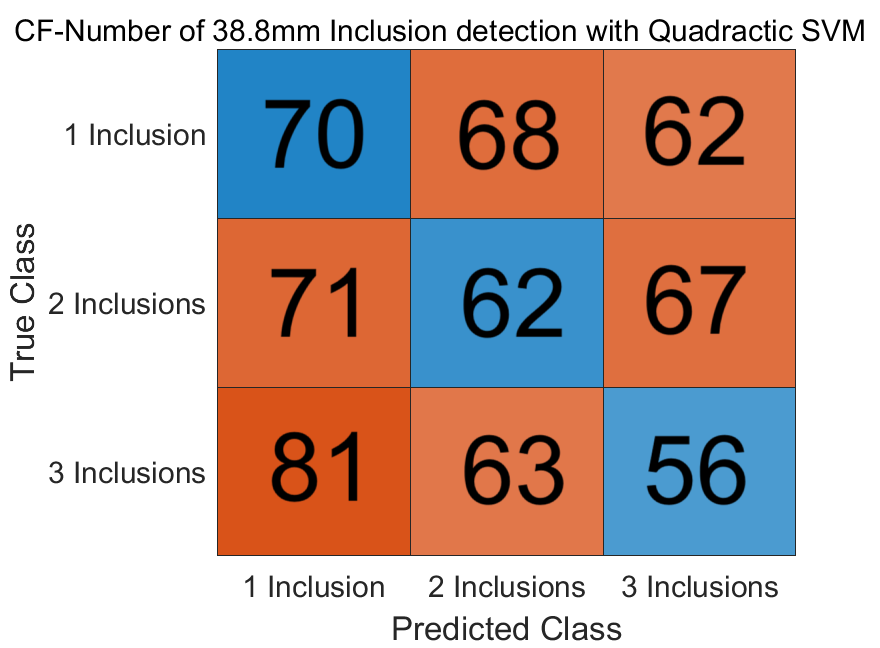}
    \caption{CF for detection the number of inclusions of radius $=38.8$mm showing $31.3\%$ accuracy in prediction.}
    \label{fig:NumInc388}
\end{figure}

Here, $600$ D-N maps (out of the $6000$ of the simulated data set) are used for testing and a $5$-fold validation is carried out as well, achieving an accuracy of $32.2\%$ on the $5$-fold validation dataset. When applied to the test dataset, the trained SVM yields an accuracy of $32.2\%$. The confusion matrix for the test data is given in Figure \ref{fig:QuadSVM_NumInc}.
This result is significantly worse than that of the SVM on the single inclusion detection problem of Section \ref{sec:1 incl}, which achieved $72.9\%$. This discrepancy indicates that the SVM is unable to determine the number of inclusions in this scenario. While the presence of a single inclusion creates a detectable difference in the D-N matrix that can be identified by the SVM, the difference caused by the presence of more inclusions (up to three small inclusions in this experiment) is not detectable by this means. 
This difference in performance between the single inclusion detection problem and the problem of determining the number of inclusions in $\Omega$ may be attributed to several factors: 
 \begin{enumerate}
    \item Since detecting the exact number of inclusions is inherently more complex than identifying their mere presence (there are more classes here), it might require a more advanced SVM, such as one with a Radial Basis Function (RBF) kernel, which could be able to capture non-linear relationships that the quadratic SVM might have missed;
    \item we tested the performance accuracy of training the inclusion number detection dataset using an ANN, as described in Section \ref{subsec:ann}, and achieved an accuracy of 34.0\%. This suggests that a more complex model, such as a Convolutional Neural Network (CNN), may be necessary;
    \item The issue could be that the D-N matrix does not vary significantly for small inclusions, thus making it challenging for the model to distinguish between different inclusion counts. 
\end{enumerate}
To address point 3, we repeat the experiment using larger inclusions with a radius of $38.8$mm (Class 2 from Table \ref{tb:class}) in experiment 2. 
We generate $6000$ samples, evenly distributed across cases with 1, 2, and 3 inclusions ($2000$ samples for each case, as before). The data is trained using the quadratic SVM, consistent with previous experiments and the classifications are as in schema~\eqref{eq:NumInc}. The resulting confusion matrix for this setup is shown in Figure \ref{fig:NumInc388}, showing a performance of $31.3\%$, similar to that in the first experiment of this section, with smaller inclusions of a fixed radius $19.4$ mm. 
Therefore, the more likely cause of the low performance is the inappropriate choice of algorithm (point 2. above) and further investigation is needed to explore alternative algorithms and approaches to ameliorate the use of the D-N matrix in ML algorithms. For example, more complex models, such as CNNs, are known to outperform simpler algorithms like SVM in classification tasks. However, we cannot entirely rule out the possibility that a different type of SVM, such as one using a cubic or RBF  kernel, might yield better performance (point 1. above). 


\section{Inclusion Radii Detection}\label{sec:IncRad}
In this section, we make the following assumption.
\begin{assumption}\label{assumption 4}
There is exactly one inclusion $D$ in $\Omega$ and it is circular.
\end{assumption}
Here we use the extensive real dataset provided by \cite{data}, comprising a circular saline tank containing a circular inclusion, to detect the radius of the inclusion via the ANN of Section \ref{subsec:ann}. 

The dataset includes a circular inclusion with four possible radii, as listed in Table \ref{tb:class}. Given that the tank consists of a saline solution, which is a uniform and homogeneous fluid, and the inclusion is composed of a single material, we assume the conductivity of both the tank and the inclusion are isotropic and constant of types \eqref{eq:itank} and \eqref{eq:iinc}, respectively. The dataset \cite{data} adopted here, contains a total of 38300 samples across all classes and is based on the so-called opposite injection pattern. We do not add noise to this data, as it is not simulated; the noise is inherently included in the measurements gathered. 

\begin{figure}[H]
    \centering
    \includegraphics[width=1\linewidth]{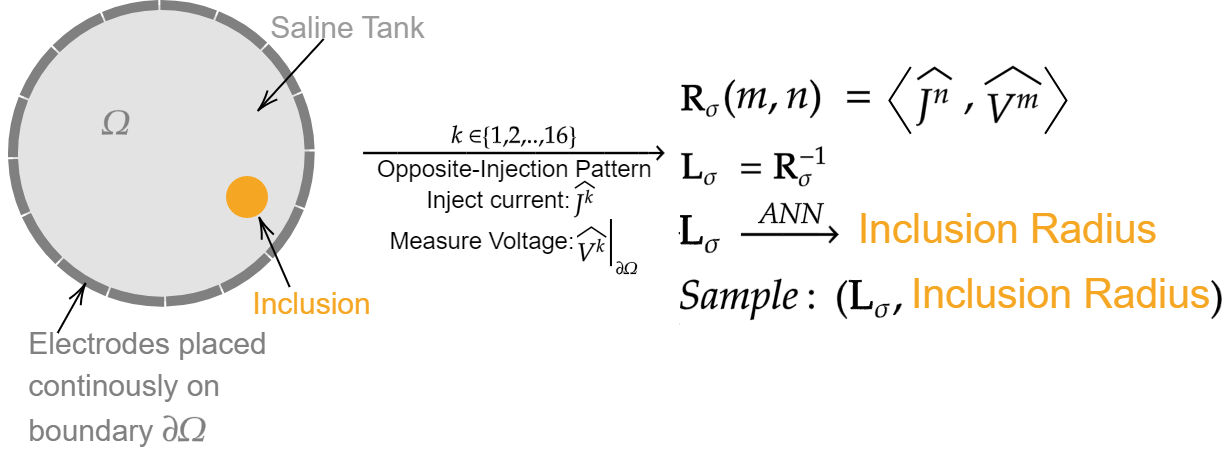}
    \caption{This figure depicts the tank and the inclusion present in the dataset. Here we use the opposite injecting pattern. Utilizing voltage and current data pairs, we generate the D-N matrix $\mathbf{L}_\sigma$, which we subsequently map to inclusion radii using our ANN. Our ANN is tasked with detecting whether the D-N matrix corresponds to a specific inclusion radius.}
    \label{fig:prob}
\end{figure}
Each $\mathbf{L_{\sigma}}$, constructed via \eqref{eq:D-N}, corresponds to a conductivity of type $\sigma = \sigma^{iso}_{inc}\:\chi_{D} + \sigma^{iso}_{tank}\:\chi_{\Omega\setminus\overline{D}}$ (here $\chi_A$ represents the characteristic function of set $A$), where the circular inclusion $D$ has radius equal one of the possible values listed in Table \ref{tb:class}, leading to samples of type $\left(\mathbf{L}_\sigma, \text{Inclusion Radius}\right)$. The D-N matrix is then classified as:
 \begin{equation}
\mathbf{L_{\sigma}} \quad\xrightarrow[]{ANN} \quad
    \begin{cases}
    \text{Inclusion diameter = 19.4mm} &  \text{Class 1}\\
    \text{Inclusion diameter = 38.8mm} &\text{Class 2}\\
       \text{Inclusion diameter = 58.2mm} & \text{Class 3} \\
       \text{Inclusion diameter = 77.6mm} & \text{Class 4}
\end{cases}.
\label{eq:radii_detection_classes}
\end{equation}
\begin{algorithm}[H]
\caption{Computation of D-N Matrix for Real Data}
\begin{algorithmic}[1]
\FOR{each sample}
    \FOR{$k = 1$ to $16$}
        \STATE For opposite current injection $J^k$ tabulate corresponding induced voltages $V^k$.
    \ENDFOR
    \STATE Use $J^k$ and $V^k$ to compute the N-D Matrix: $\mathbf{R}_\sigma$.
    \STATE Compute D-N matrix as $\mathbf{L}_\sigma=\mathbf{R}^{-1}_\sigma$
    \STATE Do not add noise to real data.
\ENDFOR
\STATE Continue until $\mathbf{L}_\sigma$ matrix for all samples is found.
\STATE Each data point consists of \{$\mathbf{L}_\sigma$, radius of inclusion in setup\}.
\STATE Train ANN using data to classify  $\mathbf{L_{\sigma}} \quad\xrightarrow[]{\text{ANN}} \text{Inclusion Radius}$.
\end{algorithmic} \label{alg:d-n-real}
\end{algorithm}

\begin{figure}[htb]
    \centering
\includegraphics[width=0.75\linewidth]{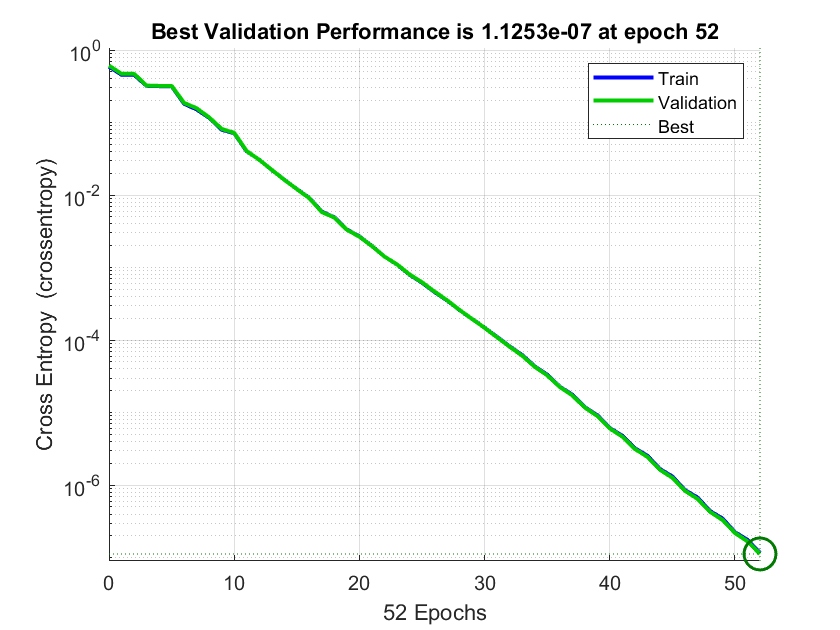}
    \caption{Cross-entropy error per epoch for training, validation and test data. We observe that the best performance on the validation dataset is reached at epoch 52. The cross-entropy error on the validation dataset is 1.1253e-07 at this point. The considerable overlap between the validation and training performance curves indicates a well-generalised model with minimal overfitting. The similarity of these two curves suggests that the model is learning patterns that generalise well to unseen data. The green circle indicates the point of optimal performance.}
    \label{fig:ecross}
\end{figure}

\begin{figure}[htb]
    \centering
\includegraphics[width=0.35\linewidth,keepaspectratio]{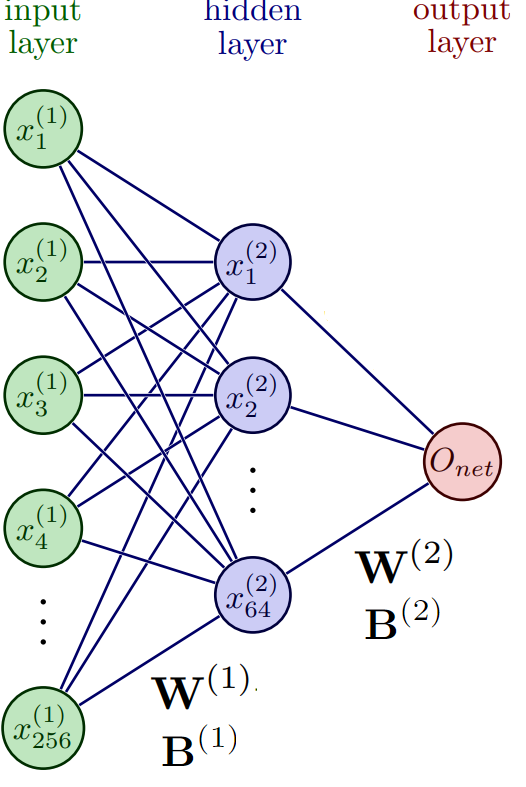}
    \caption{Architecture for R-FCNN (Radius-Fully Connected Neural Network). Neurons are represented in circles and the lines represent weights and biases. We have an input layer with 256 neurons and a hidden layer with 64 neurons followed by an output layer with 1 neuron. $\mathbf{W}^{(1)} \in \mathbb{R}^{256\times 64}$, $\mathbf{B}^{(1)} \in \mathbb{R}^{1\times64}$ and $\mathbf{x}^{(1)} \in \mathbb{R}^{1\times256}$  are the weights, biases and inputs associated with the first layer. 
     $\mathbf{W}^{(2)} \in \mathbb{R}^{64\times n}$,  $\mathbf{B}^{(2)} \in \mathbb{R}^{1\times n}$ are the weights and biases for the second layer for $n$ neurons in the output layer. The input to the second layer is defined in Equation \eqref{eq:input_nn}.
    Here $\phi$ is the sigmoid activation function in \eqref{eq:sig}, and $\varphi$ is the softmax function in  \eqref{eq:soft}. $O_{net}$ is the predicted output of the network.}
    \label{fig:RFCNN}
\end{figure}

We chose 6000 data points from each class in Table \ref{tb:class}, resulting in a total of 24000 samples in our unbiased dataset. While this is smaller than the initial dataset, it is still adequate for our ANN with one hidden layer. This is because our architecture is a simple feed-forward structure which requires far less data than more complicated networks like CNNs. 
To create a standard split for training, testing, and validation, we allocate 80\% of the dataset (19200 samples) for training, 10\% of the dataset  (2400 samples) for validation and the last 10\%  (2400 samples) of the dataset for testing as per Remark \ref{datasplitSVM}. We use the validation test to determine the end of training as shown in Figure \ref{fig:ecross}. The ANN used here is a Radius-Fully Connected Neural Network (R-FCNN), see Figure \ref{fig:RFCNN}), with the cross-entropy loss function introduced in \eqref{eq:loss}. We initially designated a maximum of 1000 epochs (iterations) for training as it is the default in MATLAB2023b, yet our algorithm reaches convergence at 52 epochs (Figure \ref{fig:ecross}), terminating the training process. When the validation error (the error in the validation dataset) consistently decreases, we continue training until reaching the maximum specified number of epochs. Since we reached a steady validation error, we completed training.

The significant improvement in results compared to those obtained in Section \ref{sec:num} could be the result of the presence of more data. Larger datasets enhance the learning abilities of ML algorithms, making distinctions clearer and improving the model's accuracy.  It will be part of future work to check whether the difference in performance might be due to the D-N matrices themselves, by varying perhaps more significantly for different inclusion radii as opposed to multiple inclusions of the same radii.

\begin{figure}[htb]
    \centering
    \includegraphics[width=0.5\linewidth]{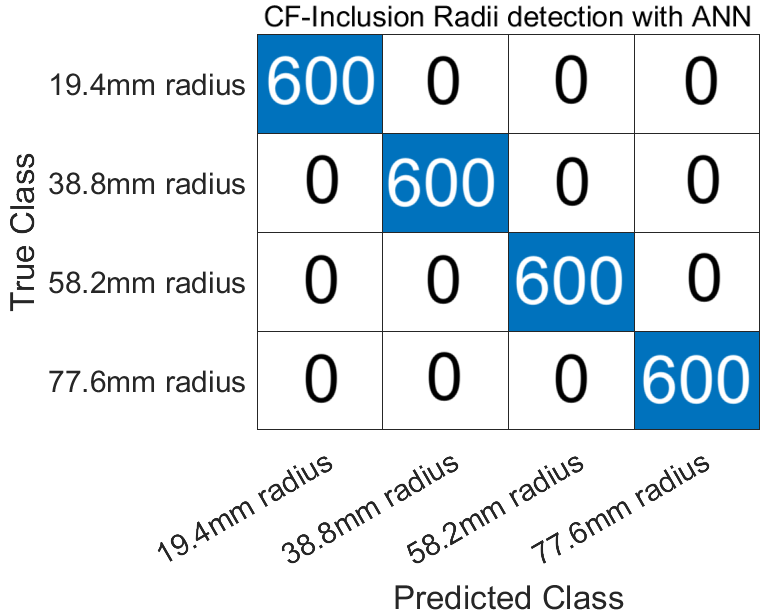}
    \caption{Confusion matrix for test and validation data showing 100\% accuracy on test data for radii's detection.}
    \label{fig:TestC}
\end{figure}

\begin{remark}
The use of real data mitigates the risk of neural networks overfitting to simulated data and prevents them from learning artifacts or biases introduced due to machine simulations, ensuring more robustness and generalisability of the trained neural networks. This emphasis on real data not only enhances the validity of our findings but also strengthens the applicability of our methodology to real-world scenarios, e.g., medical imaging applications or geophysical prospection. 
\end{remark}
\begin{remark}
We emphasise that training is performed with each sample having the inclusion positioned at a different location. By varying the location of the inclusion in each sample, the neural network is trained to accurately determine the radius of the inclusion, regardless of its position within the sample. 
\end{remark}
\begin{remark}
We employed the ANN of Section \ref{subsec:ann} because the problem is more complex as we have more classification labels, and neural networks offer strong generalisation capabilities, especially in the presence of a large dataset. 
\end{remark}
\begin{remark}
Achieving an equal distribution of samples across each class was prioritised over using all the available data for the neural network to prevent the emergence of bias towards any specific class. This emphasis on equal representation is fundamental in ML classification, including:
\begin{enumerate}
\item imbalanced classes can lead the model to favour the majority class, impeding its ability to accurately discern minority classes;
\item neural networks trained on imbalanced data may struggle to generalise effectively, particularly when the class distribution in the test set differs from that in the training set. 
\end{enumerate}
\end{remark}
By ensuring balanced representation, the network can better learn features that are relevant across all classes, thereby fostering improved generalisation \cite{Bu}. We reserved 10\% of the dataset for validation and another 10\% for testing as per Remark \ref{datasplitANN}. When tested on the held-out test dataset the ANN achieves an accuracy of 100\%, as shown in Figure \ref{fig:TestC}. We have a validation accuracy of 97.2\%. The similarity in accuracy between the testing and validation data indicates a well-generalised model, showcasing consistent performance across different datasets.

\begin{remark}
    The experiment of classifying the D-N matrix $\mathbf{L}_{\sigma}$ in the inclusion's radius experiment of this section, has also been performed by feeding an R-FCNN network with the N-D matrix $\mathbf{R}_\sigma$ directly, avoiding inverting $\mathbf{R}_\sigma$ to form  $\mathbf{L}_{\sigma}$. This was driven by the fact that the data used in this section only is real, derived by an experiment that provides directly $\mathbf{R}_\sigma$. The decision to start with inverting $\mathbf{R}_\sigma$ to form $\mathbf{L}_{\sigma}$ was driven by the sake of consistency with the remaining experiments of this manuscript aiming at classifying $\mathbf{L}_{\sigma}$ rather than $\mathbf{R}_\sigma$. When employing the N-D matrix (instead of the D-N matrix) and R-FCNN, we achieved a 100\% accuracy as well, as expected. 
\end{remark}

Next, since from a mathematical point of view, the size of an inclusion can be estimated by one measurement only (see \cite{Al-Ro-Seo}, \cite{Br}, \cite{Fr}), we test the accuracy of the inclusion radii prediction by varying both the number of measurements and electrodes in our setup. Specifically, we investigate the following two cases:
\begin{enumerate}
    \item reducing the number of measurements $M$ taken from our 16-electrode setup, from $M=6$ to $M=1$;\label{enum:var-meas}
    \item reducing the number of electrodes $E$ used to take our measurements from $E=16$ to $E=2$.\label{enum:var-elecs}
\end{enumerate}

\begin{figure}[H]
    \centering
    \includegraphics[width=0.6\linewidth]{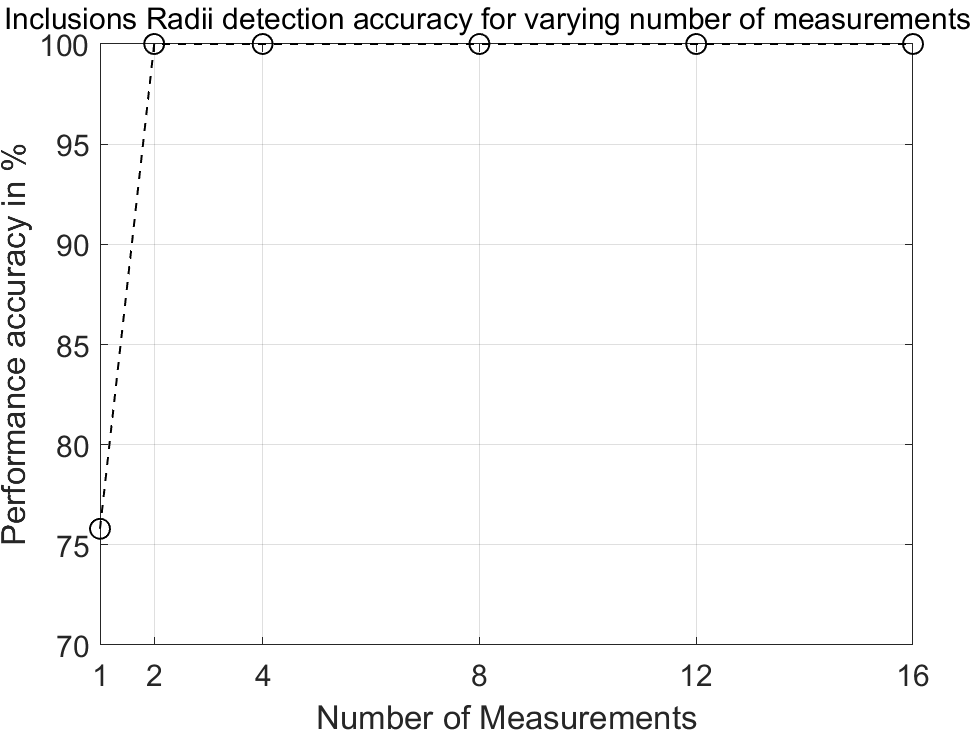}
    \caption{Performance accuracy for inclusion radii detection when varying the number of measurements $M$ taken from a 16 electrode setup and $M=1,\: 2,\: 4,\:8,\:12,\:16$.}
    \label{fig:multiple_measure}
\end{figure}

In scenario \ref{enum:var-meas}, we utilise the real data provided in \cite{data} and form the D-N matrix as previously described using \eqref{eq:D-N}. Here, with a number of measurements $M$, the D-N matrix $\mathbf{L}_{\sigma}\in \mathbb{R}^{M \times M}$, for $M=1,\: 2,\: 4,\:8,\:12,\:16$. The D-N matrix is reshaped into a column vector to obtain $\mathbf{L_\sigma} \in \mathbb{R}^{M^2}$ and the ANN, as described in Section \ref{subsec:ann}, is employed to classify $\mathbf{L_\sigma} \in \mathbb{R}^{M^2}$ as in \eqref{eq:radii_detection_classes}. As the number of neurons in the input layer of our ANN architecture is now $M$, we have $\mathbf{x}^{(1)} \in \mathbb{R}^{M^2}$, for $M=1,\: 2,\: 4,\:8,\:12,\:16$. The performance accuracy of our classifications in terms of the number of measurements $M$, $M=1,\: 2,\: 4,\:8,\:12,\:16$, from a 16-electrode system is shown in Figure \ref{fig:multiple_measure}.
Figure \ref{fig:multiple_measure} clearly shows that the use of fewer measurements does not significantly affect the results until we reach the case of a single measurement. With just one measurement, the prediction accuracy for inclusion radii drops to 75.8\%. This indicates that while a single measurement can reliably identify the inclusion radius most of the time, at least two measurements are required to achieve 100\% accuracy in identifying the inclusion radii via an ANN. Next, we consider scenario \ref{enum:var-elecs} where we decrease the number of electrodes in our setup and consider cases $E=2,\:4,\:8,\:12,\:16$.
\begin{figure}[H]
    \centering
    \includegraphics[width=0.6\linewidth]{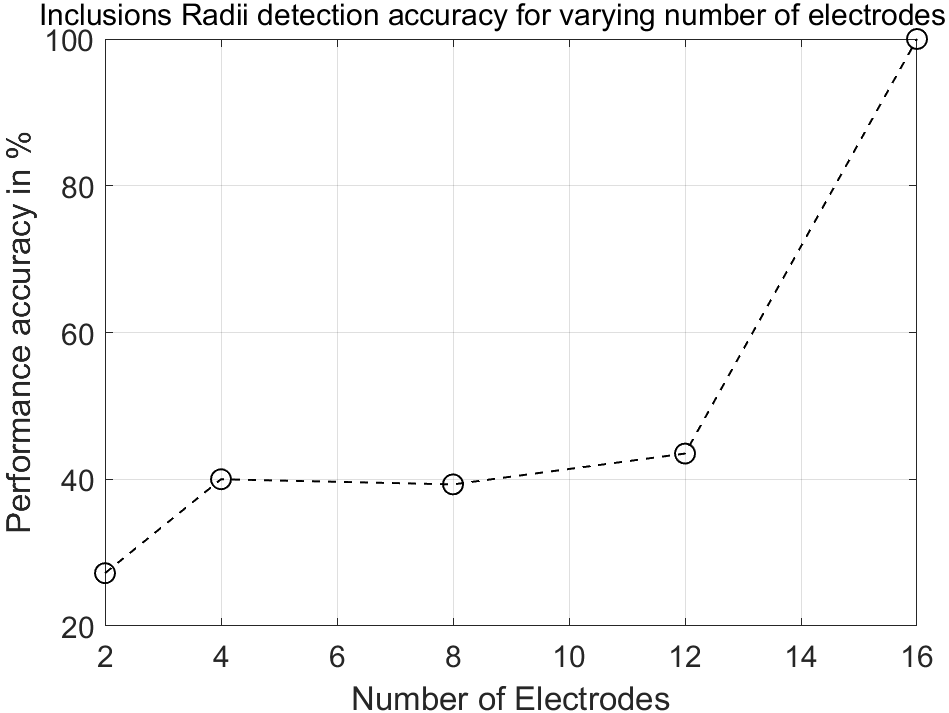}
    \caption{Performance accuracy of inclusion radii detection with varying number of electrodes. Here we take $E$ measurements from an $E$ electrode system yielding $\mathbf{L'}_\sigma \in \mathbb{R}^{E^2}$, for $E=2,\:4,\:8,\:12,\:16$.}
    \label{fig:all_elec}
\end{figure}
Since experimental data for $E<16$ are unavailable, we simulate the data for a saline tank containing an inclusion of one of the four radii listed in Table \ref{tb:class}. For each configuration with $E$ electrodes placed around the boundary of the tank, we apply $E$ opposite voltages along the boundary and compute the corresponding current densities as in section \ref{subsec:electrodes}. Unique opposite-injection patterns are not possible for $M>E$, therefore $M\leq E$, yielding to $\mathbf{L}_\sigma \in \mathbb{R}^{E\times E}$, as in \eqref{eq:D-N}. We simulate $12000$ data points for each electrode configuration, i.e., for $E=2,\:4,\:8,\:12$. Next, $\mathbf{L}_\sigma$ is reshaped into a column vector $\in \mathbb{R}^{E^2}$. We then add random noise to produce \eqref{eq:ND-N}, resulting in $\mathbf{L}'_\sigma \in \mathbb{R}^{E^2}$ which is used as input data to train the ANN as described in Section \ref{subsec:ann}, with input layer $\mathbf{x}^{(1)} \in \mathbb{R}^{E^2}$, for $E=2,\:4,\:8,\:12$. The network classifies the D-N matrix as per table \eqref{eq:radii_detection_classes}. The resulting classification accuracy displayed in Figure \ref{fig:all_elec}, clearly shows that when using a decreasing number of electrodes, the performance accuracy reduces considerably for fewer electrodes.

\begin{figure}[H]
    \centering
    \includegraphics[width=0.6\linewidth]{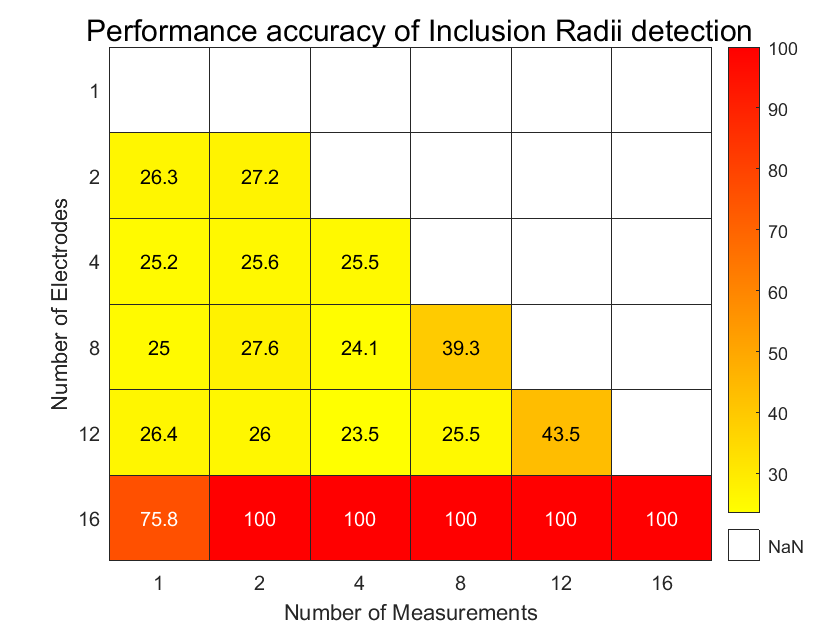}
    \caption{Heatmap of performance accuracy for Inclusion Radii detection. Setups with a single electrode are clearly not feasible, leaving the top row empty. Similarly, having the number of measurements $M$ exceed the number of electrodes $E$ is not possible, as this would result in non-unique opposite injections. Consequently, the elements above the diagonal remain empty.}
    \label{fig:heatmap}
\end{figure}

By comparing Figure \ref{fig:multiple_measure} and Figure \ref{fig:all_elec}, we observe a significant discrepancy in performance as the number of measurements and electrodes decreases, respectively. Although a 16-electrode set-up is able to achieve the precision of 100\% with only two measurements at hand, the accuracy drops to 27.2\% when using 2 measurements from a 2-electrode set-up, showing that achieving higher accuracy with fewer measurements requires a greater number of electrodes. A heatmap of performance accuracy for the number of electrodes vs the number of measurements is given in Figure \ref{fig:heatmap}.


\section{Anisotropy detection}
\label{sec:AnisoEIT}
To determine the feasibility of ANNs in detecting the presence of anisotropy within an inclusion in $\Omega$, we perform several experiments. We \textit{a-priori} assume that the saline tank $\Omega$ contains an inclusion that can be of any one radius displayed in Table \ref{tb:class}. 


\subsection{Constant diagonal anisotropy}\label{diag aniso}
Here, we define the \textit{isotropic samples} to be those samples simulating an isotropic tank containing an isotropic inclusion having conductivities \eqref{eq:iinc} and \eqref{eq:itank}, respectively; similarly, the \textit{anisotropic samples} are those samples simulating an anisotropic tank containing an anisotropic inclusion having conductivities \eqref{eq:anisomat} and \eqref{eq:atank}, respectively. To maintain unbiased samples, we generated an equal number of samples for both anisotropic and isotropic cases. We create 1000 isotropic samples, with each radius from Table \ref{tb:class} equally represented, resulting in 250 samples for each radius. Similarly, we generate 250 samples for each radius in the anisotropic case (see Table \ref{tb:allclass}). 


\subsubsection{Constant isotropic samples vs diagonal anisotropic samples}\label{classification 1}
\begin{figure}[H]
    \centering
    \includegraphics[width=0.9\linewidth]{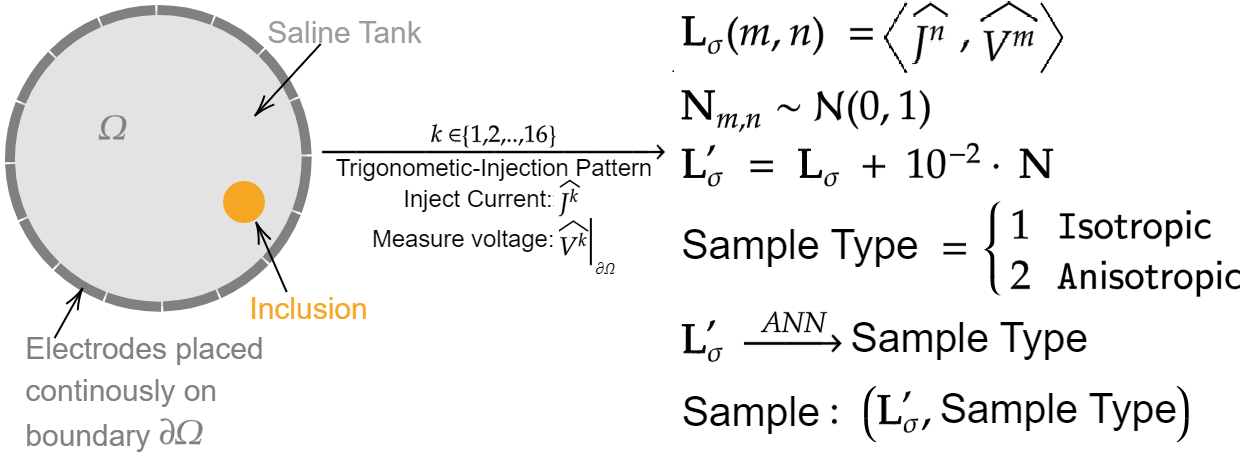}
    \caption{Schematic representation of the process followed for anisotropy detection. This is an adapted version of Figure \ref{fig:prob} to show the differences in the setup of the two problems. }
    \label{fig:prob2}
\end{figure}

We attribute Class `1' for isotropic samples and for anisotropic samples, we attribute Class `2'. The ANN of Section \ref{subsec:ann} is tasked with determining which class the noisy D-N Matrix $\mathbf{L'}_\sigma$ in \eqref{eq:ND-N}
maps to, i.e. 
$$
    \mathbf{L}^{'}_{\sigma} \xrightarrow{\text{R-FCNN}} \begin{cases}
        \text{Isotropic sample} & \text{Class 1} \\
        \text{Anisotropic sample}& \text{Class 2}
    \end{cases}
    \label{eq:goal2}.
$$
We use the trigonometric bases voltage of \eqref{eq:vol}. The D-N matrix is calculated as in \eqref{eq:D-N} and noise is added to the D-N matrix to get $\mathbf{L'}_\sigma$ as in \eqref{eq:ND-N}, see Figure \ref{fig:prob2}.
\begin{remark}
Simple ANNs (like R-FCNN) have fewer parameters than more complex architectures like CNNs or Recurrent Neural Networks, making them less data-hungry and suitable for smaller datasets \cite{cnn}. This feature of ANNs makes them ideal for the training of 2000 samples only in this section.
\end{remark}
The isotropic constant conductivity of the saline tank and inclusions is \eqref{eq:itank} and \eqref{eq:iinc}, with $\lambda_{tank}=1.45 \,\text{S} \,\text{m}^{-1}$ \cite{saline}, while the constant conductivity 
of the inclusions is set at $\lambda_{inc}=10 \, \text{S} \,\text{m}^{-1}$ based on the exhibited conductivity of thermoplastics used in 3D printing, particularly those incorporating carbon nanotubes \cite{pla}.
For the anisotropic experiment, the anisotropic tank and inclusion is given by 
\begin{equation}
 \sigma^{aniso}_{tank}  = \mu_{tank}
\begin{bmatrix}
a & 0 \\
0 & b  \\
\end{bmatrix},
\label{eq:atank}
\end{equation}
and
\begin{equation}
 \sigma^{aniso}_{inc}  = \mu_{inc}
\begin{bmatrix}
a & 0 \\
0 & b  \\
\end{bmatrix},
\label{eq:anisomat}
\end{equation}
respectively, where $2 \leq a,b \leq 10 $ are uniformly distributed pseudorandom integers generated using \texttt{randi(2,10)} in MATLAB2023b. This was done to ensure $\sigma^{aniso}$ is a symmetric positive definite matrix \eqref{ellipticity anisotropic}. We simulate our data using  \eqref{eq:vol} and the continuum electrode model of Section \ref{subsec:electrodes}. We calculate the noisy D-N Matrix $\mathbf{L}^{'}_{\sigma}$ as per \eqref{eq:D-N}-\eqref{eq:ND-N}. 
 \begin{figure}[H]
    \centering
    \includegraphics[width=0.5\linewidth]{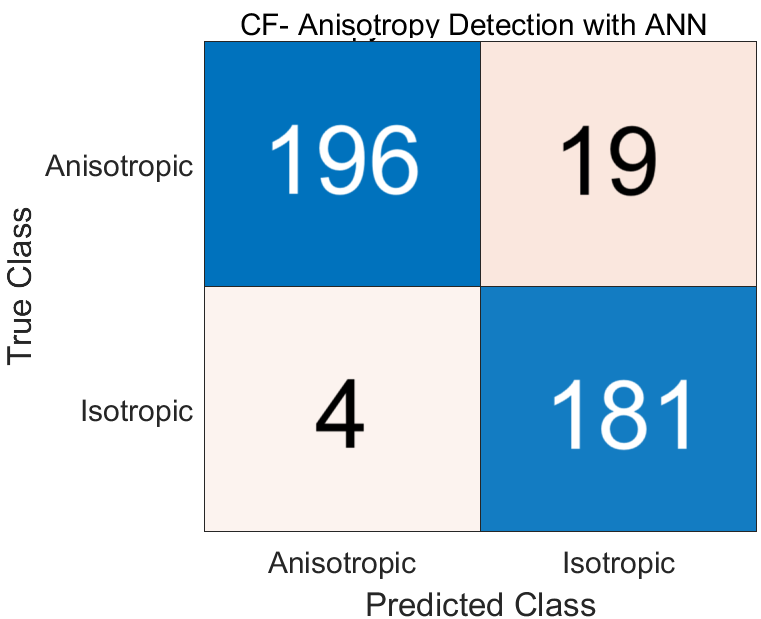}
    \caption{Confusion matrix for anisotropy detection. Class 1 refers to an isotropic tank and inclusion, and Class 2 refers to an anisotropic tank and inclusion. The observed testing accuracy stands at $94.2\%$, denoting a $94.2\%$ detection rate of anisotropy. }
    \label{fig:anisodet}
\end{figure}
As for the inclusion radii detection classification problem of Section \ref{sec:IncRad}, we employ the R-FCNN of Section \ref{subsec:ann}. This is to facilitate a performance comparison between inclusion radii detection and inclusion anisotropy detection considered here. The only difference between these two classification problem is the input data and the number of output classes. We had 4 classes for inclusion radii detection, each for one radius. Here we have 2 classes, one for isotropic samples and one for anisotropic samples, as defined above. We initialise the weights and biases arbitrarily and train the network on the simulated dataset. The results of training R-FCNN on this new dataset are shown in Figure \ref{fig:anisodet}. Anisotropy is here predicted with an accuracy of $94.2\%$, compared to $100.0\%$ accuracy for inclusion radii detection. We observe a test accuracy of $94.2\%$. Additionally, it's reassuring to note a close prediction accuracy of $96.7\%$ on the validation dataset. This suggests that our data is not overly prone to overfitting.


\subsubsection{Constant isotropic inclusion vs anisotropic inclusion}\label{classification 2}
\begin{figure}[H]
    \centering
    \includegraphics[width=0.5\linewidth]{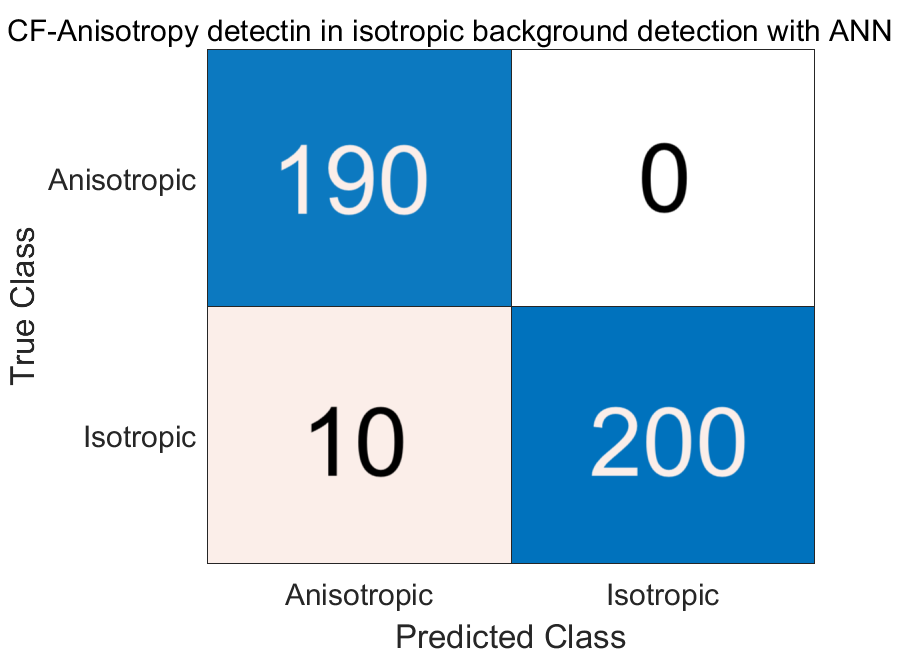}
    \caption{CF for detection of anisotropic inclusion in an isotropic tank.}
    \label{fig:IsoTankANisoIso}
\end{figure}

Next, we extend our analysis to determine whether an inclusion within an isotropic tank is isotropic or anisotropic. For this experiment, we use an isotropic tank of the form \eqref{eq:itank} containing an inclusion, where the inclusion could be either isotropic of form \eqref{eq:iinc} or anisotropic of form \eqref{eq:anisomat}. The ANN's task is to detect the presence of anisotropy in the inclusion, i.e.
$$
\mathbf{L}^{'}_{\sigma} \xrightarrow{\text{R-FCNN}} \begin{cases}
        \text{Isotropic inclusion} &\text{Class 1} \\
        \text{Anisotropic inclusion} & \text{Class 2}
    \end{cases}
    \label{eq:goal2b}.
$$
    The resulting confusion matrix for this scenario is shown in Figure \ref{fig:IsoTankANisoIso}, showing that the presence of an anisotropic in an isotropic background inclusion can be detected with an accuracy of 97.5\%.


\subsection{Constant diagonal vs off-diagonal anisotropy}\label{off-diag aniso}
Next, we consider anisotropic conductivities that are both diagonal and off-diagonal types of conductivities, replicating the setup of the real dataset in \cite{data}. Specifically, we generate 4000 samples exhibiting diagonal anisotropy. In these diagonal anisotropic samples, both the tank and the inclusion are diagonally anisotropic as in \eqref{eq:atank},\eqref{eq:anisomat}. We also generate an additional 4000 samples that we call \textit{off-diagonal anisotropy}. For these samples, both the tank and the inclusion exhibit off-diagonal anisotropy as in \eqref{eq:atankoff} and \eqref{eq:anisomatOffDiag} below, respectively.
\begin{equation}
 \sigma^{aniso}_{tank}  = \mu_{tank}
\begin{bmatrix}
a & c \\
c & b  \\
\end{bmatrix},
\label{eq:atankoff}
\end{equation}
and
\begin{equation}
 \sigma^{aniso}_{inc}  = \mu_{inc}
\begin{bmatrix}
a & c \\
c & b  \\
\end{bmatrix},
\label{eq:anisomatOffDiag}
\end{equation}
respectively, where $6 \leq a,b \leq 20 $ and $1\leq c\leq5$ are uniformly distributed pseudorandom integers generated using \texttt{randi()} in MATLAB2023b to ensure $\sigma$ is a positive-definite matrix.  The inclusion, as before, is characterized by one of the four radii described in Table \ref{tb:class}. We apply the trigonometric voltage pattern as outlined in \eqref{eq:vol}, calculate the $\mathbf{L}_{\sigma}$ D-N matrix as defined in \eqref{eq:D-N}, and then add noise to produce the noisy D-N matrix $\mathbf{L}^{'}_{\sigma}$ as in \eqref{eq:ND-N}. As with previous setups, we allocate 10\% of the data for validation and another 10\% for testing. Finally, we train the ANN, as detailed in Section \ref{subsec:ann}, to classify the input D-N matrix as: 
$$
    \mathbf{L}^{'}_{\sigma} \xrightarrow{\text{R-FCNN}} \begin{cases}
        \text{Diagonal anisotropy} & \text{Class 1} \\
        \text{Off-diagonal anisotropy} & \text{Class 2}
    \end{cases}
    \label{eq:goal3}.
$$

\begin{figure}[htb]
    \centering
    \includegraphics[width=0.5\linewidth]{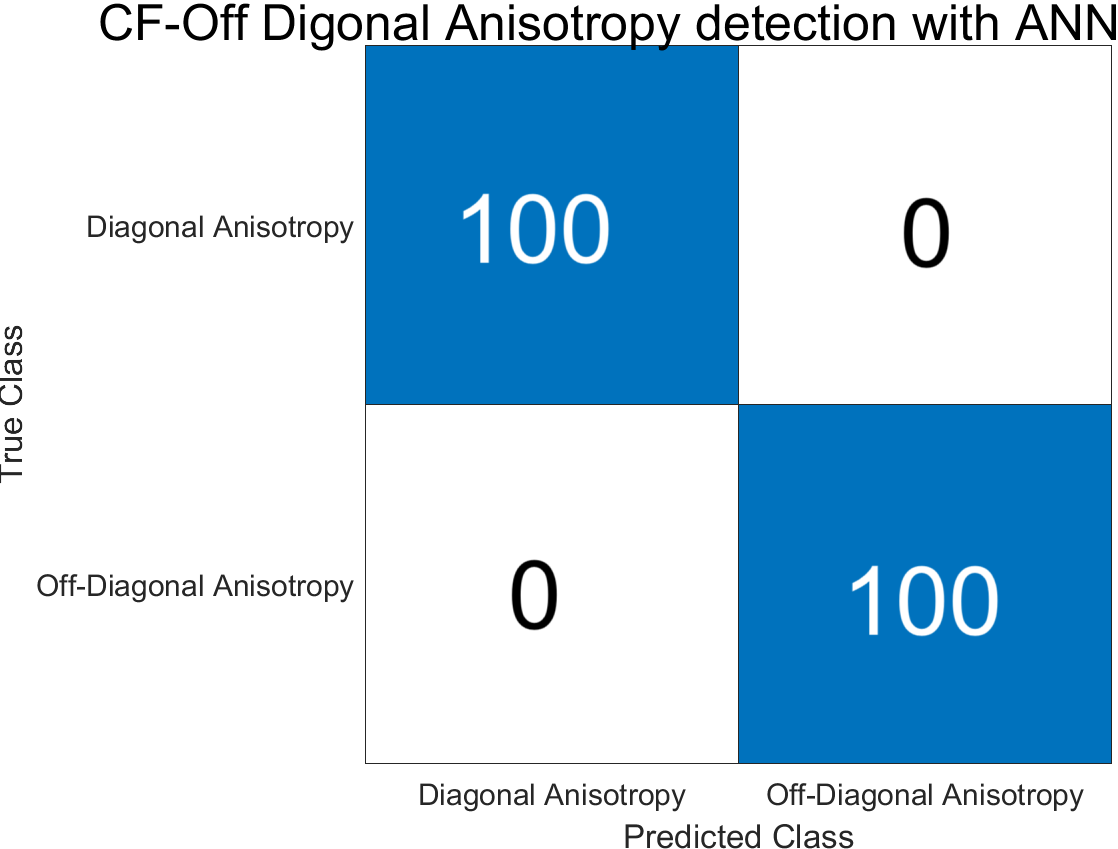}
    \caption{CF for Off-Diagonal Anisotropy detection using ANN.}
    \label{fig:OffDiag}
\end{figure}

Using R-FCNN the accuracy achieved here is 100\%. The CF for this is given in Figure \ref{fig:OffDiag}.

\subsection{Constant isotropic inclusion vs spatially varying anisotropic inclusion}\label{sec: last}
In this last experiment, we investigate the ANN's performance in classifying the D-N matrix corresponding to either of the cases: 
\begin{enumerate}
\item isotropic samples as per Subsection \ref{diag aniso}, where the tank is isotropic with conductivity \eqref{eq:itank} and the inclusion is isotropic with conductivity \eqref{eq:iinc};\\
\item anisotropic spatially varying samples where the tank is isotropic and constant of type \eqref{eq:itank} and the inclusion is anisotropic, spatially varying, of type
\begin{equation}
 \sigma^{aniso}(x,y)  = 
\begin{bmatrix}
x^2 & 0 \\
0 & y^2  \\
\end{bmatrix}.
\label{eq:spatialcond}
\end{equation}
\end{enumerate}
Data is simulated as in the rest of Section \ref{sec:AnisoEIT}, placing at most one circular inclusion within the tank $\Omega$, having one of four possible radii as in Table \ref{tb:class}. The D-N matrix is simulated as in Subsections \ref{diag aniso} and \ref{off-diag aniso} above employing 16 electrodes, yielding to $\mathbf{L}_\sigma \in \mathbb{R}^{16 \times 16}$ and reshaped in a vector of $\mathbb{R}^{256} $. Noise is added to get $\mathbf{L'}_\sigma \in \mathbb{R}^{256}$ as in \eqref{eq:ND-N}. A total of 8000 samples is generated. Of these, 4000 samples feature isotropic samples. The remaining 4000 exhibited spatially varying conductivity for the inclusion of the form \eqref{eq:spatialcond} with an isotropic saline tank of the form \eqref{eq:itank}.

We trained the R-FCNN as described in Section \ref{subsec:ann} to classify $\mathbf{L'}_\sigma$ as,
$$
    \mathbf{L}^{'}_{\sigma} \xrightarrow{\text{R-FCNN}} \begin{cases}
        \text{Isotropic inclusion} & \text{Class 1} \\
        \text{Anisotropic spatially varying inclusion} & \text{Class 2}
    \end{cases}
    \label{eq:goalspatial}.
$$
\begin{figure}[htb]
    \centering
    \includegraphics[width=0.5\linewidth]{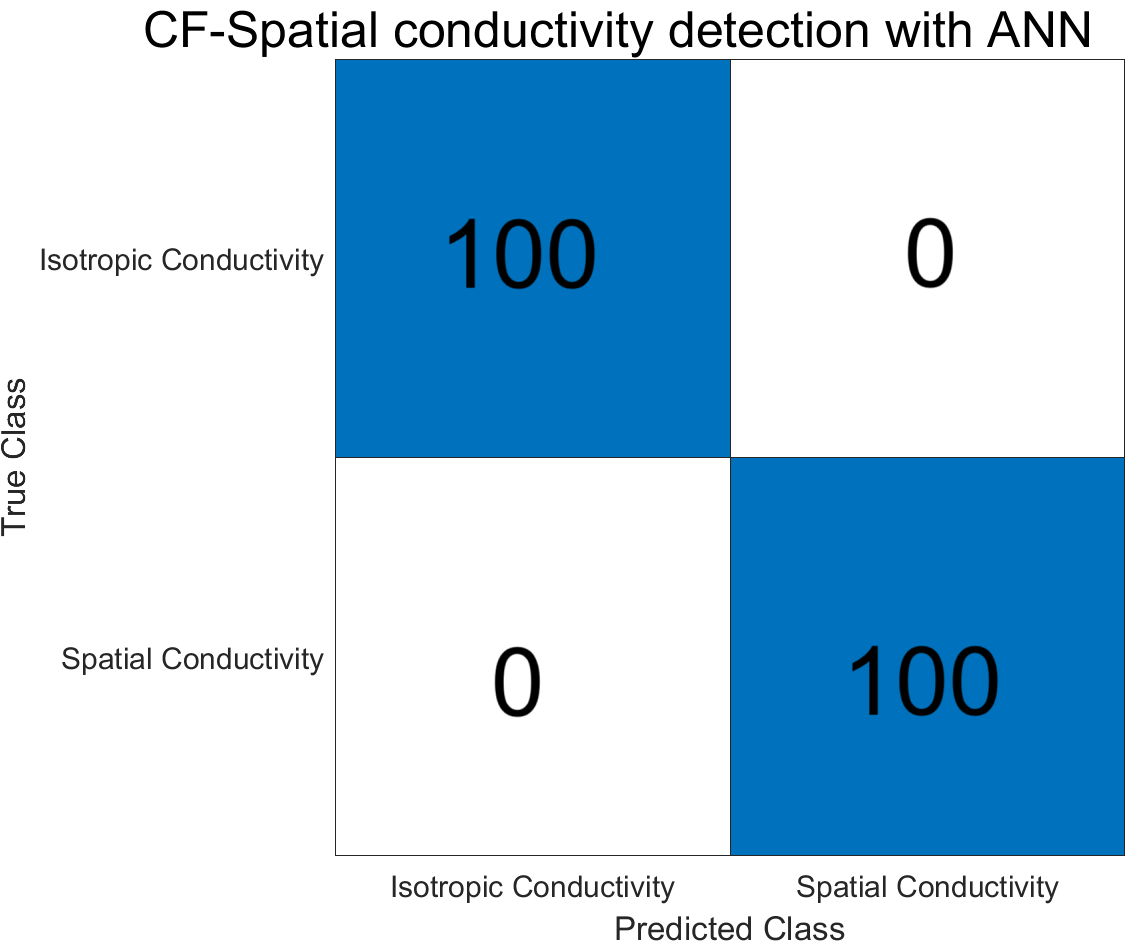}
    \caption{CF for spatially varying conductivity detection}
    \label{fig:SpatialCond}
\end{figure}
Consistently with the rest of the manuscript, we use 80\% of the data for training. We use 10\% of the remaining data for validation and the last 10\% for testing. We achieve an accuracy of 100.0\% in the detection of spatially varying conductivity. The CF for this is given in Figure \ref{fig:SpatialCond}.

To summarize this section, we build upon the findings of \cite{dC-R1} by considering different conductivities for the inclusion $ D $ and the background $ \Omega \setminus \overline{D} $. While \cite{dC-R1} focuses on the case where the inclusion exhibits a conductivity jump across $ \partial D $ relative to the background, our study introduces four additional distinct types.  The specific differences in conductivity we examine are outlined in Table \ref{tb:AnisoType}.

\begin{table}[H]
\centering
\begin{tabular}{lrr}
\textbf{Anisotropy Type}                       & \textbf{Anisotropic Inclusion} $D$          & \textbf{Background} $\Omega\setminus\overline{D}$ \\ \hline
Constant diagonal                  & $\sigma=diag(a,b) $ \eqref{eq:anisomat}        &  $\sigma=diag(a,b) $ \eqref{eq:atank}                \\
Constant diagonal                  & $\sigma=diag(a,b) $ \eqref{eq:anisomat}       & $\sigma=\lambda \mathbf{I}$ \eqref{eq:itank}                 \\
Off-Diagonal &  $\sigma=sym(a,c,b) $
 \eqref{eq:anisomatOffDiag} &$\sigma=sym(a,c,b) $
 \eqref{eq:atankoff}             \\
Spatially Varying         & $\sigma=diag(x^2,y^2) $ 
 \eqref{eq:spatialcond}     & $\sigma=\lambda \mathbf{I}$ \eqref{eq:itank}               
\end{tabular}
\caption{Table of types of anisotropic inclusions and the respective background considered in this manuscript, where $sym()$ is a symmetric square matrix. }
\label{tb:AnisoType}
\end{table}


\section{Conclusions}
\label{sec:con}
The application of ML techniques to the inverse problem of detecting an inclusion and its conductivity in a conducting body $\Omega\subset\mathbb{R}^2$ from electrostatic measurements, is explored. Specifically, using both real and simulated datasets of the D-N matrix corresponding to a $16$ electrodes setup around $\Omega$, the following problems are addressed: i) detecting the presence of an inclusion $D$ within $\Omega$ (Section \ref{sec:1 incl}); ii) determining the number of inclusions (Section \ref{sec:num}); iii) identifying inclusions radii (Section \ref{sec:IncRad}); iv) detecting the presence of anisotropy within $D$ (Section \ref{sec:AnisoEIT}). 

Two ML methods are employed: SVMs (Section \ref{subsec:SVM_alg}) and ANNs (Section \ref{subsec:ann}), allowing for ML classification tasks in i)-iv) above. The results are summarised in Table \ref{tb:results}.

\begin{table}[H]
\centering
\begin{tabular}{lll}

\textbf{D-N matrix and ML algorithm's task}                                             & \textbf{Method} & \textbf{Result} \\ \hline
One inclusion $D$ detection                         & SVM   & 72.9\%            \\
Inclusions number (19.4mm)                & SVM   & 32.2\%            \\
Inclusions number (38.8mm)                & SVM   & 31.3\%            \\
Radius of $D$                             & ANN     & 100.0\%           \\
Constant isotropy vs diagonal anisotropy  & ANN      & 94.2\%            \\
Constant isotropic $D$ vs anisotropic $D$  & ANN       & 97.5\% \\
Constant diagonal vs off-diagonal anisotropy & ANN             & 100.0\% \\
Constant isotropic $D$ vs spatially varying anisotropic $D$ & ANN & 100.0\% \\
\end{tabular}
\caption{Summary of results in this manuscript}
\label{tb:results}
\end{table}

It is anticipated that this proof of concept paper will serve as a first step for future work on identifying anisotropic conductivities within a body. This will include expanding our work to the treatment of more general inhomogeneous (spatially varying) anisotropic conductivities, where the conductivity of the inclusion $D$ differs from that of the background $\Omega$ not only by a (constant) scalar factor. This would allow to fully extend the theoretical results obtained in \cite{dC-R1} via ML techniques. Here, a partial extension of \cite{dC-R1} is obtained (see Table \ref{tb:AnisoType}). Specifically, in Section \ref{sec: last}, the classification task involves the distinction between a spatially varying anisotropic conductor $D$ from a constant isotropic one in a constant isotropic conducting background $\Omega$. 

As the current paper is concerned with the 2D problem only, another interesting direction will be to explore the 3D case as well. These extensions will possibly be explored via more advanced networks like CNNs to handle added complexity. 

Improvements in SVM performance, through optimised kernels or alternative methods is another promising direction the authors wish to explore as part of their future work to improve the low accuracy obtained here when classifying the number of inclusions within $\Omega$ (see Section \ref{sec:num}). 

Another interesting problem would be that of using more realistic electrode models in EIT, such as the complete electrode model, together with ML techniques like decision trees to broaden the ML potential in EIT data recognition. 
\renewcommand{\thefootnote}{}
\footnotetext{The code written and the data simulated for this manuscript is made available at: https://zenodo.org/records/14778868.}


\section*{Acknowledgments}
This publication has emanated from research conducted with the financial support of Science Foundation Ireland under Grant number 18/CRT/6049. For the purpose of Open Access, the author has applied a CC BY public copyright license to any Author Accepted Manuscript version arising from this submission. 
\\
Neural Network diagrams (Figures \ref{fig:an}, \ref{fig:RFCNN}) were made using \cite{tikz_neural_networks}.

\end{document}